\documentclass[11pt]{amsart}
\usepackage{amsmath,amssymb,amsthm,amscd}

\def\p{\partial}
\def\R{\mathbb{R}}

\def\i{\sqrt{-1}}
\def\l{\lambda}

\def\D{\Delta}

\def\cE{{\mathcal E}}

\def\cR{{\mathcal R}}
\def\cS{{\mathcal S}}

\def\tr{\text{Tr}}

\newtheorem{prop}{Proposition}[section]
\newtheorem{thm}[prop]{Theorem}

\newtheorem{rmk}[prop]{Remark}
\newtheorem{exam}[prop]{Example}

\newtheorem{pl}[prop]{Problem}

\begin{document}
\title{Nonlinear Hodge flows in symplectic geometry}
\author{Weiyong He}

\begin{abstract}
Given a symplectic class $[\omega]$ on a four torus $T^4$ (or a $K3$ surface), a folklore problem in symplectic geometry is whether symplectic forms in $[\omega]$ are isotropic to each other. We introduce a family of nonlinear Hodge heat flows on compact symplectic four manifolds to approach this problem, which is an adaption of nonlinear Hodge theory in symplectic geometry. 
As a particular example, we study a conformal Hodge heat flow in detail. We prove a stability result of the flow near an almost Kahler structure $(M, \omega, g)$. We also prove that, if $|\nabla \log u|$ stays bounded along the flow, then the flow exists for all time for any initial symplectic form $\rho\in [\omega]$ and it converges to $\omega$ smoothly along the flow with uniform control, where $u$ is the volume potential of $\rho$. 
\end{abstract}

\maketitle
\section{Introduction}
 Let $(M, \omega)$ be a compact symplectic four manifold such that $a=[\omega]\in H^2(M; \R)$. Denote $\cS_a$ to be the set of the symplectic forms in the class of $a$,
\[
\cS_a=\left\{\rho\in \Omega^2(M): \rho^2>0, d\rho=0, [\rho]=a\right\}. 
\]
A long standing open problem in symplectic topology is 
\begin{pl}\label{connected}Is $\cS_a$ path connected?
\end{pl}
The uniqueness question for the four torus $T^4$ is a long standing question since 1980s. 
It is generally believed that Problem \ref{connected} holds true for  hyperK\"ahler surfaces (and maybe $\mathbb{CP}^2$ as well). A related problem on K\"ahler surfaces was proposed by S. Donaldson \cite{Donaldson06},
\begin{pl}\label{q1}Suppose $X$ is a compact K\"ahler surface with K\"ahler form $\omega_0$. If $\omega$ is any other symplectic form on $X$, with the same Chern class and with $[\omega] = [\omega_0]$, is there a diffeomorphism f of X with $f^*\omega=\omega_0$?
\end{pl}
If Problem \ref{connected} is confirmed then Problem \ref{q1} follows as a consequence via the Moser isotropy. 
It is known, by the deep theory in symplectic topology, that certain symplectic structures in the same cohomology class in dimension four are unique modulo diffeomorphisms, including the projective surface $\mathbb{CP}^2$, see \cite[Section 13.4]{MS} for a detailed discussion of examples and references. Hence in particular if Problem \ref{connected} is confirmed for $\mathbb{CP}^2$, it would imply that every diffeomorphism of $\mathbb{CP}^2$ that acts as the identity on homology is isotopic to the identity. However Problem \ref{connected} remains completely open for any compact symplectic four manifolds.
We refer to the survey paper by  D. Salamon \cite{Salamon} and the book \cite[Chapter 13, 14]{MS} by D. McDuff and D. Salamon for an excellent overview. 

S. Donaldson described a geometric flow approach for Problem \ref{connected}  in \cite{Donaldson99}. We roughly recall the setup and refer the readers to \cite{KS} for a detailed discussion. Associate $(M, \omega)$ with a compatible Riemannian metric $g$. For any $\rho\in \cS_a$, define
\[
\text{dvol}_\rho:=\frac{\rho\wedge \rho}{2}, u:=\frac{\rho\wedge \rho}{\omega\wedge \omega}=\frac{\text{dvol}_\rho}{\text{dvol}_g}, g(A \cdot, \cdot):=\rho, g_\rho:=u^{-1} g(A\cdot, A\cdot)
\]
The metric $g_\rho$ induces the Hodge star $*_\rho$ and the Donaldson metric on the tangent space at $\rho$, $T_\rho \cS_a$. The Donaldson flow, defined by
\begin{equation}\label{Donaldsonflow}
\p_t \rho=d*_\rho d \left(\text{II}(u^{-1} *\rho)\right), \;\text{II}(\tau)=\tau-\frac{\tau \wedge \rho}{\rho\wedge \rho} \rho,\; \text{for}\; \tau \in \Omega^2(M).  
\end{equation} can be viewed as the negative gradient flow of the energy functional (with respect to the Donaldson metric), 
\[
\cE(\rho)=\int_M \frac{|\rho|^2}{u} dv,
\] 
Note that  $u^{-1}|\rho|^2\geq 2$ pointwise and it achieves the absolute minimum when $*\rho=\rho$. Hence $\omega$ is the unique minimum of $\cE(\rho)$ over $a=[\omega]$. 
If one could prove that the Donaldson flow exists for all time and converges to $\omega$ for a \emph{generic} $\rho\in a$, Problem \ref{connected} can be answered affirmatively.  
For the progress in the Donaldson's geometric flow, see \cite{KS, Krom}.

In this paper we introduce a family of parabolic flow of symplectic forms, motivated by the Hodge theory, to approach Problem \ref{connected}. Let $(M, \omega, g)$ be an almost K\"ahler manifold of dimension four. Given a symplectic form $\rho\in [\omega]$, one can associate a Riemannian metric $\tilde g$, depending on $\rho$ and $g$. In this paper we mainly discuss a conformal metric $\tilde g=\sqrt{u} g$, where $u$ is the volume potential of $\rho$ such that
\[
u=\frac{\rho\wedge \rho}{\omega\wedge \omega}=\frac{\text{dvol}_\rho}{\text{dvol}_g}, \text{with}\; \text{dvol}_\rho=\frac{\rho\wedge \rho}{2}
\]
Denote $ *$ to be the Hodge star of $g$, and $\tilde *$ the Hodge star of $\tilde g$. Given a choice of $\tilde g=\tilde g(\rho, g)$, 
we consider the parabolic flow of $\rho$ by
\begin{equation}\label{general-flow}
\p_t \rho=d\tilde * d*\rho.
\end{equation}
There is a flexibility to choose the metric $\tilde g$. 
For any such choice, the energy functional 
\[
E(\rho)=\int_M \rho\wedge *\rho
\]
is decreasing along the flow and $\omega$ is the unique minimizer of $E(\rho)$, by the Hodge theory. In dimension four, \eqref{general-flow} has a description as a gradient flow of $E(\rho)$ if one uses the metric $\tilde g$ to define the metric on the tangent space of $\rho\in \cS_a$. \\

When $\tilde g=g_u:=\sqrt{u} g$, the flow equation reads
\begin{equation}\label{conformal-flow}
\p_t \rho=-d\left(\frac{d^*\rho}{\sqrt{u}}\right)
\end{equation}

The first result regarding  \eqref{conformal-flow} is the following stability near $(M, \omega, g)$, 
\begin{thm}\label{main-1}There exists $\epsilon>0$ depending on $(M, \omega, g)$, such that for any $\rho_0\in [\omega]$ satisfying  $\|\rho_0-\omega\|_{C^{1, \alpha}}<\epsilon$, then \eqref{conformal-flow} exists for all time with initial data $\rho_0$ and converges to $\omega$, at the exponential decay rate when $t\rightarrow \infty$. 
\end{thm}

\begin{rmk}With reasonable effort, Theorem \ref{main-1} should be true by assuming $\rho_0\in [\omega]$ such that for some $p>4$, 
$
\|\rho_0-\omega\|_{W^{1, p}}<\epsilon. $
It remains a challenge problem to prove the stability if we assume $\rho_0\in L^\infty\cap W^{1, 4}$ such that
$\|\rho_0-\omega\|_{W^{1, 4}}<\epsilon. $
\end{rmk}

One key point for Theorem \ref{main-1} is to prove that the normalized energy functional
\[
E_0(\rho)=E(\rho)-E(\omega)
\]
is decreasing along the flow at exponential rate, see Proposition \ref{energydecay}. The second result is the following long time existence and convergence criterion, 
\begin{thm}\label{main-2}Suppose $|\nabla \log u|$ stays uniformly bounded along the flow, then \eqref{conformal-flow} has a long time solution such that $\rho$ converges to $\omega$ smoothly with uniform control along the flow. The convergence  is at exponential decay rate, for a set of suitable positive constants $A_0, A_1, C$, satisfying the estimate,
\begin{equation}\label{edecay-01}
\left(\int_{M} |d^*\rho|^2+A_1 |\rho_0-\omega|^2\right)(t) \leq A_0 e^{-Ct}. 
\end{equation}
\end{thm}

With the uniform control on $|\rho|$, $|\nabla \rho|$ and exponential decay \eqref{edecay-01}, it is a standard practice to get higher order exponential decay on $\|(d^*d)^k \rho\|_{L^2}$ for each $k\geq 1$. 
The key to Theorem \ref{main-2} is to explore the nonlinear structure of the equation and consider the function $f$ of the form, for suitable positive constants $A$ and $B$, 
\begin{equation}\label{convergence-01}
f=|\nabla \rho|^2+A|\nabla u|^2+B |\rho|^2+1.
\end{equation}
Assuming the uniform control of $|\nabla \log u|$, one can obtain a differential inequality for $f$,
\begin{equation}\label{convergence-02}
\p_t f\leq \frac{\Delta f}{\sqrt{u}}+C_1 f.
\end{equation}
A maximum principle argument then implies that $f(t)\leq C_0 e^{C_1 t}$. 
With this bound, it is now standard to obtain a Shi-type estimate on higher derivatives of $\rho$, and obtain the long time existence. The exponential decay estimate \eqref{edecay-01} is a refinement from the exponential decay of energy. 
To obtain the smooth convergence along the flow with uniform control, we need to bound $|\rho|$ and $|\nabla \rho|$ pointwise. The key point is to use the differential inequality \eqref{convergence-02} to bound $f$, via a delicate parabolic Moser iteration. Hence if $|\nabla \log u|$ is bounded along the flow, Theorem \ref{main-2} would solve Problem \ref{connected} with affirmative answer. 
However, this remains a very challenge problem. Instead of assuming $C^1$ bound on the potential $u$, it would be a very interesting question to relate the long time existence of the flow to the geometry of the conformal metric $g_u=\sqrt{u} g$. 

It is not clear  whether \eqref{conformal-flow} will develop a (finite time) singularity or not in dimension four. Our result implies that $|\nabla \log u|$ blows up at a possible finite time singularity. 
If the dimension of $M$ $n=6$, D. McDuff \cite{McDuff} constructed an example on $T^2\times S^2\times S^2$ such that there is a symplectic cohomology class, that $\cS_a$ is not connected. One can consider a geometric flow, for any $r>0$
\begin{equation}\label{rflow}
\p_t \rho=-d\left(\frac{d^*\rho}{u^r}\right).
\end{equation}
The energy functional $E(\rho)$ is still decreasing. If \eqref{rflow} exists for all time with uniform control, then the flow will converge to $\omega$, the unique minimizer of $E(\rho)$. This implies that \eqref{rflow} would have to develop singularities for McDuff's example. 
One might expect that there would be singularities along the flow \eqref{conformal-flow}  in general. If that were the case, it would be an important question to study the formation of singularity of \eqref{conformal-flow},
\begin{pl}Does the flow \eqref{conformal-flow} form a (finite time) singularity or the flow exists for all time in general?
\end{pl}

The nonlinear structure of the conformal Hodge flow is of the form of a nonlinear harmonic heat flow type, with a conformal factor $u^{-\frac{1}{2}}$ on the leading term, 
\[
\p_t \rho=\frac{\Delta \rho}{\sqrt{u}}+\nabla (u^{-\frac{1}{2}})\circ \nabla \rho+\text{l.o.t}
\]
There are obviously two major difficulties in the study of \eqref{general-flow} and \eqref{conformal-flow}. 
The first one  is whether the conformal factor $u^{-\frac{1}{2}}$ stays uniformly positive and bounded along the flow. 
The second one is related to the nonlinear structure of the harmonic heat flow type. These two difficulties are intervened, and it remains a very interesting and challenge problem to go beyond Theorem \ref{main-2}, with more flexible assumptions on $u$.

One can compare with the linear Hodge flow in this setting,
\[
\p_t \rho=-d(d^*\rho)=d*d*\rho
\]
By the standard parabolic theory and the Hodge theory, the flow exists for all time, with uniform control on $\rho$, and converges to $\omega$, the unique minimizer of $E(\rho)$. The potential function $u$ stays uniformly bounded and converges to $1$. However, $u$ might not stay uniformly positive and hence $\rho$ would not stay as a nondegenerate two form along the linear Hodge flow. Hence the linear Hodge flow is not a right tool to approach Problem \ref{connected}. Our method is a nonlinear modification of the linear Hodge flow.  One might ask, whether these is a choice of $\tilde g=\tilde g(\rho, g)$, such that $u$ stays uniformly positive along the flow \eqref{general-flow}. 

\begin{pl}Is there a particular choice of $\tilde g$, such that \eqref{general-flow} exists for all time and converges to $\omega$?
\end{pl}

We will discuss some specific choices of $\tilde g$ in Section \ref{example2}.  Our computations suggest it is a very delicate problem. In particular, we construct an example (Example \ref{negative-01}) of a sequence of symplectic forms on $T^4$ with the normalized volume potential $1$, along the linear Hodge flow, the volume potential can turn negative almost instantly (the time can be as small as possible for the sequence). \\

{\bf Acknowledgement:} It is a pleasure to acknowledge the valuable discussions with Song Sun over the years, in particular regarding the problems raised in \cite{Donaldson06}. 

\numberwithin{equation}{section}

\section{A nonlinear Hodge theory in symplectic geometry}

We introduce a family of geometric flows for Problem \ref{connected}, motivated by the Hodge theory. Such flows can also be defined in higher dimensions. 
Let $(M, \omega, g)$ be a compact symplectic manifold with a compatible metric $g$ of dimension $2n$.
By the Hodge theory,  $\omega$ is the unique harmonic form in the class $a=[\omega]$, which minimizes the norm over $a$, 
\[
E(\alpha)=\int_M |\alpha|^2 dv=\int_M \alpha\wedge *\alpha
\]
The motivation is then to minimize $E(\rho)$ over $\cS_\alpha$ while enforcing the nondegenerate condition $\rho^n>0$. Suppose $\delta \rho= d\xi,$ for $\xi\in \Omega^1(M)$. The variation of $E(\rho)$ is given by
\[
\delta E(\rho)=\int_M \xi \wedge d*\rho.
\]
Let $\tilde g_\rho$ be a Riemannian metric induced by $g$ and $\rho$, with the Hodge star $\tilde *$. We compute
\[
\delta E(\rho)=-\int_M \xi \wedge \tilde * (\tilde * d*\rho)=-\int_M \langle \xi, \tilde * d*\rho\rangle d\tilde v. 
\]
By taking $\xi=\tilde * d*\rho$ and hence $\delta \rho=d(\tilde * d*\rho)$, we obtain
\[
\delta E(\rho)=-\int_M |\tilde * d*\rho|^2_{\tilde g_\rho} d\tilde v.
\]
We consider the following flow
\begin{equation}\label{gradient1}
\p_t \rho=d\left(\tilde * d*\rho\right).
\end{equation}
If we fix $\tilde g_\rho=g$, then \eqref{gradient1} is reduced to the Hodge Laplacian heat flow. The Hodge Laplacian heat flow exists for all time and converges to $\omega$ for any initial two form $\alpha \in a$, but the nondegenerate condition $\rho^n>0$ cannot be assured.  
We use the choice of $\tilde g_{\rho}$ to enforce $\rho^n>0$, which is evolved along the flow \eqref{gradient1} as well. 
By choosing the metric $\tilde g_\rho$, 
\eqref{gradient1} can be viewed as a gradient flow of $E(\rho)$ if we use $g_{\tilde \rho}$ to define a metric on the forms in the tangent space at $\rho\in \cS_a$, similar as  in the setting of the Donaldson geometric flow \cite{KS}. 
Allowing to choose  $\tilde g_\rho$ seems to be valuable.  A critical point remains the same for $E(\rho)$  (or for $\cE(\rho)$). In particular, the flow \eqref{gradient1} decreases the functional $E(\rho)$ and the only fixed point of the flow satisfies $d*\rho=0$; by the Hodge theory such a fixed point $\rho$ must be a harmonic form  and equals $\omega$. Hence if the gradient flow \eqref{gradient1} has a long time solution with $C^{1, \alpha}$ norm uniformly bounded, the flow will converge to $\omega$. 

Given a symplectic form $\rho$ and a metric $g$, there are many ways to produce a (family of) Riemannian metric $\tilde g_\rho$. We define the linear transformation $A: TM\rightarrow TM$ by
\[
\rho=g(A\cdot, \cdot)
\]
It follows that $B=-A^2$ is symmetric and positive definite, and hence $B^r$ is uniquely well-defined for any real number $r$. Note that $B^r A=AB^r$.  Suppose the eigenvalues of $A$ reads $\{\pm \i \l_1, \cdots, \pm \i \l_n\}$ with $\l_i>0$ for $i=1, \cdots, n$. The eigenvalues of $B$ read $\{\l_1^2, \cdots, \l_n^2\}$, with each repeated. 
We define  $\tilde g_{\rho}=g(Q \cdot, \cdot)$ with $Q=Q(\l_1^2, \cdots, \l_n^2)$ to be a positive definite symmetric matrix such that $QA=AQ$. A typical choice of $Q$ is of the form $fB^r$ for a positive function $f=f(\l_1, \cdots, \l_n)$ and $r\in \R$.  Since our primary interest is on symplectic manifolds of dimension four, we focus on the case $n=2$. 

\begin{prop}\label{shorttime}Let $(M, \omega, g)$ be an almost Kahler manifold of dimension four. For any initial data $\rho_0\in a$, \eqref{gradient1} has  a unique smooth short time solution $\rho(t)$ in $[0, t_0)$ with $\rho(0)=\rho_0$, for any choice of $Q$. 
\end{prop}
If we apply Hamilton-Nash-Moser implicit function theorem \cite{Hamilton}, the short time existence of \eqref{gradient1} follows immediately since the system has an obvious integrability operator $L=d$ and $d\tilde *d*$ is an elliptic operator on $\text{Ker}(d)$. Here we give an elementary argument for the short time existence and for simplicity, we consider only $n=2$. 
Consider at the moment $\rho$ is only a nondegenerate two form. We consider
\begin{equation}\label{parabolic1}
\p_t \rho=d\tilde * d*\rho+  *d\tilde * d\rho
\end{equation}

\begin{prop}Let $(M, g)$ be a four manifold with a nondegenerate two form  $\rho_0$ such that $\rho_0^2>0$. Then for and any choice of $\tilde g_\rho$ described above, the flow \eqref{parabolic1} is strictly parabolic and admits a unique short time smooth solution. 
Along the flow we have
\begin{equation}\label{decreasing}
\begin{split}
\p_t \int_M \rho\wedge *\rho=&-\int_M d\rho \wedge \tilde * d\rho-\int_M (d*\rho)\wedge \tilde * (d*\rho)\\
=&-\int_M \left(|d\rho|^2_{\tilde g}+|d*\rho|^2_{\tilde g}\right) d\tilde v
\end{split}
\end{equation}
Moreover if $d\rho_0=0$, then $d\rho=0$ along the flow. 
\end{prop}

\begin{proof} 
Recall $\tilde g_\rho=g (Q\cdot, \cdot )$ with $Q=Q(\l_1^2, \l_2^2)$. For any $\theta\in \Omega^3(M)$, we compute
\[
\tilde *\theta=(\det{Q})^{-\frac{1}{2}} (*\theta) \circ Q. 
\]
We have the following,
\[
\tilde *\circ *|_{\Lambda^1}=*\circ \tilde *|_{\Lambda^3}=-(\det{Q})^{-\frac{1}{2}} Q
\]
Denote $\tilde Q=(\det{Q})^{-\frac{1}{2}} Q.$ We rewrite \eqref{parabolic1} as
\[
\p_t \rho=- d(\tilde Q \circ d^*\rho)-d^*(\tilde Q \circ d\rho)
\]
 We only need to consider the leading terms and hence at one point we can assume that $\tilde Q_{kl}=\tilde \l_k \delta_{kl}$.
Then the leading term reads
\[
\p_t \rho_{ij}=\sum_k\tilde \l_k \p^2_k \rho_{ij}+\text{l.o.t}
\]
and hence \eqref{parabolic1} is a strictly parabolic equation on $\rho$. This proves the short time existence when $n=2$. 
Hence for any initial nondegenerate two form $\rho_0$, there exists a unique smooth short time solution $\rho(t)$ in $[0, T)$. For any $t_0\in (0, T), \rho(t)$ is a smooth nondegenerate two form in $[0, t_0].$ If $d\rho_0=0$, we compute
\[
\begin{split}
\p_t \int_M |d\rho|^2 dv=&2\int_M  \langle d * d\tilde *d\rho, d\rho \rangle dv\\
=&2\int_M \langle  * d\tilde*d\rho, d^*d\rho\rangle dv\\
=&-2\int_M\langle d^* \tilde Q d\rho, d^*d\rho\rangle dv
\end{split}
\]
 Hence we compute
 \[
 \begin{split}
 \p_t  \int_M |d\rho|^2 dv=&-2\int_M \langle d^* (\tilde Q d\rho), d^*d\rho\rangle dv\\
 \leq &-c_0\int_M \langle d^*d\rho, d^*d\rho\rangle dv+C_1\int_M |d\rho||d^*d\rho|dv\\
 \leq& \;\;C_2\int_M |d\rho|^2 dv
 \end{split}
 \]
 where $c_0$, $C_1$ and $C_2$ are positive constant depending on $t_0$.  Gronwall's inequality implies that $d\rho=0$ over $[0, t_0]$. Note that along the flow \eqref{parabolic1}, the energy functional $E(\rho)$ is still decreasing such that \eqref{decreasing} holds, which is a straightforward computation. 
 
 \end{proof}

 We describe several interesting choices of $Q$, given $\rho$ and $g$. Firstly we fix notations. The almost K\"ahler structure $(\omega, g)$ is fixed. We choose an orthonormal frame $\{e_1, e_2, e_3, e_4\}$ and its dual coframe $\{\theta^1, \theta^2, \theta^3, \theta^4\}$. We write
\[
\rho=\frac{1}{2} \rho_{ij} \theta^i\wedge \theta^j=\sum_{i<j} \rho_{ij} \theta^i\wedge \theta^j
\]
Recall that $\rho=g(A\cdot, \cdot)$. Define $B$ to be $*\rho=g(B \cdot, \cdot)$
\[
A=\begin{pmatrix}
0 & \rho_{12} & \rho_{13} &\rho_{14}\\
-\rho_{12}&0 &\rho_{23}&\rho_{24}\\
-\rho_{13}&-\rho_{23} &0 &\rho_{34}\\
-\rho_{14}&-\rho_{24}&-\rho_{34}&0
\end{pmatrix},\,\,\, B=\begin{pmatrix}
0 & \rho_{34} & -\rho_{24} &\rho_{23}\\
-\rho_{34}&0 &\rho_{14}&-\rho_{13}\\
\rho_{24}&-\rho_{14} &0 &\rho_{12}\\
-\rho_{23}&\rho_{13}&-\rho_{12}&0
\end{pmatrix}
\]
We have $AB=-u\,\text{I},\;\; B=-uA^{-1}$, with $u=\rho_{12}\rho_{34}-\rho_{13}\rho_{24}+\rho_{14}\rho_{23}=\sqrt{\det{A}}$. 
Denote
\[
\rho^{+}=\frac{1}{2} \left(\rho+*\rho\right), \rho^{-}=\frac{1}{2}\left(\rho-*\rho\right)
\]
We have the following, 
\[
|\rho|^2=\sum_{i<j} \rho_{ij}^2=|\rho^+|^2+|\rho^{-}|^2, 2u=\sum_{i<j} \rho_{ij} (*\rho)_{ij}=\langle \rho, *\rho\rangle=|\rho^+|^2-|\rho^{-}|^2.
\]
\begin{prop}The eigenvalues of $A$ are $\pm\i \l_1$, $\pm \i \l_2$, with 
\[
\l_1=\frac{1}{\sqrt{2}} (|\rho^{+}|+|\rho^{-}|), \l_2=\frac{1}{\sqrt{2}} (|\rho^{+}|-|\rho^{-}|)
\]
The eigenvalues of $(a_{ij})=-A^2$ and $(b_{ij})=-B^2$ are $\{\l_1^2, \l_1^2, \l_2^2, \l_2^2\}$ and $u=\sqrt{\det{A}}=\lambda_1 \lambda_2$. Moreover,
\[
a_{ij}+b_{ij}=|\rho|^2 \delta_{ij}
\]
\end{prop}

\begin{proof}Since $A$ is skew-symmetric and nondegenerate, hence the eigenvalues of $A$ are purely imaginary,  
denoted by $\pm \i \l_1, \pm \i\l_2$. We can assume $\l_1, \l_2$ are both positive.  It is then straightforward to compute its eigenvalues in terms of $|\rho|^2$ and $u$. In particular,
\[
|\rho|^2=\l_1^2+\l_2^2,\; u=\l_1 \l_2.
\]
\end{proof}

\begin{exam}
We describe several interesting choices of $Q$.
\begin{enumerate}
\item Take $Q_{ij}=u^{-1}a_{ij}$, then $\det {Q}=1$, the metric $g_\rho=u^{-1}g(a, \cdot)=u^{-1} g(A\cdot, A\cdot)$. This is the metric used in \cite{Donaldson99} to define Donaldson's geometric flow.
One computes \cite[Proposition 3.1]{KS}
\[
*_\rho \rho=\frac{|\rho|^2}{u} \rho -*\rho,\;\text{and}\; \rho \wedge *_\rho \rho=\rho\wedge *\rho. 
\]

\item Take $Q=e^{2\varphi} \text{I}$ for a conformal factor  $\varphi=\varphi(\rho, g)$, then the metric $g_\varphi=e^{2\varphi} g$ is a conformal change. A particular interesting choice is $\varphi=\frac{\log u}{4}$, hence $\text {dvol}_\varphi=\text{dvol}_\rho$. For any conformal change, the Hodge star satisfies
\[
*_\varphi= e^{(4-2p)\varphi} *
\]
for a $p$-form. Hence
\[
\rho \wedge *_\varphi \rho=\rho\wedge *\rho. 
\]

\item There are many other choices of $Q=Q(\lambda_1, \lambda_2)$. An important question is whether there is a choice of $Q$ such that $u$ stays uniformly positive along the flow.
\end{enumerate}
\end{exam}

\section{A conformal Hodge flow in dimension four}

Let $(M, \omega, g)$ be an almost K\"ahler manifold of dimension four. 
In this section we consider $g_\varphi= e^{2\varphi} g $ with $\varphi=\frac{1}{4} \log u$, and hence $\text{dvol}_{\varphi}= u \text{dvol}_g=\text{dvol}_\rho$. In this case, since $* \rho=\tilde * \rho$, the equation \eqref{gradient1} can also be written as Hodge Laplacian flow with respect to $\tilde g$, 
\[
\p_t \rho= d\tilde * d\tilde * \rho=-\tilde{ \Delta}^d \rho
\]
We choose a normalization condition such that
\[
\text{Vol}(M, g)=\int_M \text{dvol}_g=\int_M \text{dvol}_\rho=1. 
\]
Note that the Hodge star of $e^{2\varphi} g$ is related to that of $g$ by, for a $p$ form $\theta$
\[
*_\varphi \theta = e^{(4-2p)\varphi} * \theta
\]
Hence for the conformal metric $ \sqrt{u} g$, the flow reads 
\begin{equation}\label{conformal-flow1}
\p_t \rho=-d\left(\frac{d^* \rho}{ \sqrt{u}}\right)=-\frac{dd^* \rho}{\sqrt{u}}+\frac{du\wedge d^* \rho}{2\sqrt{u^3}}=\frac{\Delta_g \rho+\cR (\rho)}{\sqrt{u}}+\frac{du\wedge d^* \rho}{2\sqrt{u^3}},
\end{equation}
where $\cR(\rho)$ is the curvature tensor acting on $\rho$. To be more precise, for an orthonormal frame $\{e_1, e_2, e_3, e_4\}$ and the corresponding coframe
$\{\theta^1, \theta^2, \theta^3, \theta^4\}$ and a form $\zeta$, 
\[
\cR(\zeta)=\theta^i\wedge \iota_{e_j} R(e_i, e_j)\zeta
\]
We use the following sign convention for the curvature operator,
\[
\begin{split}
R(e_i, e_j)e_k&=\left(\nabla_{e_i}\nabla_{e_j}-\nabla_{e_j}\nabla_{e_i}-\nabla_{[e_i, e_j]}\right) e_k\\
Ric(X, X)&=\langle R(e_i, X)X, e_i\rangle
\end{split}
\]
Denote $\D$ to be the Laplace-Beltrami operator of $g$, for any tensor $T$,
\[
\D T=\tr(\nabla^2 T)= \sum \left(\nabla_{e_i}\nabla_{e_i} T-\nabla_{\nabla_{e_i}e_i} T\right).
\]
We recall the Bochner-Weitzenb\"ock identities as follows, 
\begin{prop}\label{P-BW}
For any $k$-form $\zeta$, 
\begin{equation}\label{BW01}
\begin{split}
\D_d \zeta &=-\D \zeta-\cR(\zeta)\\
\frac{1}{2}\D \left(|\zeta|^2\right)&=\langle -\D_d \zeta, \zeta\rangle+|\nabla \zeta|^2+\langle \cR(\zeta), \zeta\rangle
\end{split}
\end{equation}
where  $\cR(\zeta)=\theta^i\wedge \iota_{e_j} R(e_i, e_j)\zeta$. 
\end{prop}

Note that  \eqref{conformal-flow1} is a quasilinear parabolic system and if  the volume potential $u$ is  uniformly positive and H\"older continuous, the system has the structure
\[
\p_t \rho= L(\rho)+\nabla \rho *\nabla \rho,
\]
which is similar to the structure of a nonlinear harmonic map heat flow, with the leading term involved with a conformal factor of $u$. A standard parabolic theory implies the following,
\begin{prop}Let $\rho_0$ be a $C^{1, \alpha}$ symplectic form in $[\omega]$. There exists a unique solution of \eqref{conformal-flow1} over a short time $[0, T)$ for some $T>0$, which is smooth for any $t\in (0, T)$. 
\end{prop}

\subsection{The stability near the energy minimizer $\omega$} In the class of $[\omega]$, $\omega$ is the only harmonic form with respect to $g$, which minimizes the energy functional. When we restrict to the symplectic forms in $[\omega]$, the flow \eqref{conformal-flow1} can be viewed as a gradient flow of the energy functional $E(\rho)$. Hence it is expected that if one starts with an initial symplectic form $\rho_0$ near $\omega$, the flow exists for all time and converges to $\omega$. We justify the expectation in this section and the key point is the following monotonicity of the normalized energy.  Denote 
\[
E_0=E(\omega)=\int_M \omega\wedge \omega=\int_M \omega\wedge *\omega=\int_M |\omega|^2 \text{dvol}_g
\]
For simplicity of the notations, we will drop the volume elements $\text{dvol}_g$ in the integrals when there exists no ambiguity. 
For a closed form  $\alpha\in [\omega]$, we have
\[
\int_M |\alpha-\omega|^2 =\int_M |\alpha|^2-2\int_M \alpha\wedge \omega+\int_M |\omega|^2=E(\alpha)-E_0
\]
Denote $E_0(\alpha)=E(\alpha)-E_0$. We have the following Poincare inequality with $c_0=c_0(g)$, 
\begin{equation}\label{poincare1}
E_0(\alpha)\leq c_0 \int_M |d^*\alpha|^2 
\end{equation}
We need the following equivalence for $W^{1, p}$ norms on the differential forms. For $p\geq 1$, one defines the geometric $W^{1, p}$ norm of a form $\alpha$ by
\[
\|\alpha\|_{W^{1, p}}^p=\int_M |\alpha|^p+\int_M |d \alpha|^p+\int_M |d^*\alpha|^p
\]
and the analytic $W^{1, p}_A$ norm is defined as
\[
\|\alpha\|_{W_{A}^{1, p}}^p=\int_M |\alpha|^p+\int_M |\nabla \alpha|^p
\]
When $p=2$, it is a standard practice that two Sobolev norms are equivalent. A key point we need is that these two norms are equivalent for any $p\geq 1$, see \cite[Corollary 12]{Scott}

\begin{prop}[C. Scott]\label{equivalent1}For any $k$-form $\alpha$ on a compact Riemannian manifold $(M, g)$, two Sobolev norms $W^{1, p}$and $W^{1, p}_A$ are equivalent for $p\geq 1$. 
\end{prop}
As a consequence of this equivalence, we have the following, 

\begin{prop}
For any close form $\alpha\in [\omega]$ and $p\geq 1$, we have the Poincare inequality, 
\begin{equation}\label{poincare2}
\int_M |\alpha-\omega|^p\leq c_1 \int_M |d^*\alpha|^p
\end{equation}
Moreover we have the Sobolev-Poincare inequality on $(M^4, \omega, g)$,  
\begin{equation}\label{sobolev1}
\int_M |\alpha-\omega|^2 \leq c_2 \left(\int_M |d^*\alpha|^{\frac{4}{3}}\right)^{\frac{3}{2}}
\end{equation}
\end{prop}

\begin{proof}
Given the equivalence of two Sobolev norms on differential forms, the proof is rather standard and we sketch a brief argument. The standard Sobolev inequality implies that
\[
\int_M |\alpha-\omega|^2 \leq C \left(\int_M |\nabla (\alpha-\omega)|^{\frac{4}{3}}\right)^{\frac{3}{2}}+C\left(\int_M |\alpha-\omega|^{\frac{4}{3}}\right)^{\frac{3}{2}}
\]
The equivalence of the Sobolev norms implies that, with $d\alpha=0$, 
\[
\int_M |\alpha-\omega|^2  \leq C\left( \int_M |d^*\alpha|^{\frac{4}{3}}\right)^{\frac{3}{2}}+C\left(\int_M |\alpha-\omega|^{\frac{4}{3}}\right)^{\frac{3}{2}}
\]
Hence \eqref{sobolev1} follows by the Poincare inequality \eqref{poincare2}, with $p=4/3$. The proof of \eqref{poincare2} follows by a standard compactness argument and norm equivalence in Proposition \ref{equivalent1}. Suppose otherwise, then there exists a sequence of forms $\alpha_k\in [\omega]$ such that 
\[
\int_M |d^*\alpha_k|^p\leq \frac{1}{k} \int_M |\alpha_k-\omega|^p
\]
Suppose $\int_M |\alpha_k-\omega|^p=\lambda_k^p$ for each $k$. We consider $\beta_k=\lambda_k^{-1} (\alpha_k-\omega)$. Then $\beta_k$ is an exact two-form  and 
\[
\int_M |d^*\beta_k|^p\leq \frac{1}{k}, \int_M |\beta_k|^p=1. 
\]
We write $\beta_k=d\zeta_k$, where  $\zeta_k\in \text{Im}(d^*)$ by Hodge decomposition.
By the equivalence of the Sobolev norms, $\beta_k\in W^{1, p}$ with bounded Sobolev norm; note that by the Hodge theory $\zeta_k$ is uniquely determined by $\beta_k$ and has uniform $W^{2, p}$ bound. Hence $\beta_k$ converges by subsequence weakly in $W^{1, p}$, to an exact form $\beta_\infty=d\zeta_\infty$.  Moreover, the weak convergence implies that 
\[
\int_M |d^*\beta_\infty|^p\leq \liminf \int_M |d^*\beta_k|^p=0.
\]
It follows that $d^*\beta_\infty=d^*d\zeta_\infty=0$ and hence $\beta_\infty=d\zeta_\infty=0$ (we have $\zeta_\infty=0$ as well, as a consequence of $\zeta_\infty\in \text{Im}(d^*)$). 
However the Rellich-Kondrachov theorem implies $\beta_k$ converges  strongly in $L^p$. 
This contradicts  the fact that $\beta_k$ converges strongly in $L^p$ to $\beta_\infty=0$, since $\|\beta_k\|_{L^p}=1$. 
\end{proof}

With the Sobolev-Poincare inequality, we have the following monotonicity of the normalized energy functional $E_0(\rho)$ along the flow,

\begin{prop}Suppose the flow \eqref{conformal-flow1} exists in $[0, T)$. For any $t<T$, 
\begin{equation}\label{energydecay}
E_0(\rho(t))\leq E_0(\rho_0) e^{-\frac{t}{c_2}}
\end{equation}
\end{prop}

\begin{proof}We have, along the flow, 
\[
\p_t \int_M |\rho-\omega|^2=\p_t \int_M |\rho|^2=-\int_M \frac{|d^* \rho|^2}{\sqrt{u}}
\]
Apply the H\"older inequality
\[
\int_M fh \leq \left(\int_M f^{\frac{3}{2}}\right)^{\frac{2}{3}} \left(\int_M h^3\right)^{\frac{1}{3}},\,\text{with}\,\,f=\frac{|d^*\rho|^{\frac{4}{3}}}{u^{\frac{1}{3}}}, h=u^{\frac{1}{3}}.
\]
We get
\[
\left(\int_M |d^*\rho|^{\frac{4}{3}}\right)^{\frac{3}{2}}\leq \int_M \frac{|d^* \rho|^2}{\sqrt{u}} \left(\int_M u\right)^{\frac{1}{2}}= \int_M \frac{|d^* \rho|^2}{\sqrt{u}}
\]
By the Sobolev-Poincare inequality \eqref{sobolev1}, 
we have 
\[
\int_M |\rho-\omega|^2\leq c_2 \left(\int_M |d^*\rho|^{\frac{4}{3}}\right)^{\frac{3}{2}}\leq c_2 \int_M \frac{|d^* \rho|^2}{\sqrt{u}}
\]
We get \eqref{energydecay} as a consequence of the following, 
\[
\p_t \int_M |\rho-\omega|^2\leq -\frac{1}{c_2} \int_M |\rho-\omega|^2
\]
\end{proof}
Now we are ready to state the main theorem in this section,
\begin{thm}[Stability] \label{stability}There exists $\epsilon>0$, such that if $\|\rho_0-\omega\|_{C^{1, \alpha}}<\epsilon$, then the flow \eqref{conformal-flow1} exists over $[0, \infty)$ and converges to $\omega$ when $t\rightarrow \infty$. 
\end{thm}

The proof of Theorem \ref{stability} follows the standard line of arguments, given the exponential decay of normalized energy in  \eqref{energydecay}. Firstly we need a more precise statement of short time existence, including the smooth dependence of initial data (finite time stability). The following is a standard fact of quasi-parabolic system applied to \eqref{conformal-flow1}. 

\begin{prop}\label{short2}There exists $\delta_0=\delta_0(\omega, g)$,  $t_0=t_0(\omega, g)$, $C_0=C_0(\omega, g)$, and $C_k=C_k(\omega, g, k)$, $k\geq 2$ for any initial data $\rho_0$ satisfying
\[
\|\rho_0-\omega\|_{C^{1, \alpha}}<\delta_0,
\]
then the flow exists over $[0, t_0]$ and for any $t\in [0, t_0]$, 
\[
\|\rho-\omega\|_{C^{1, \alpha}}\leq C_0 \|\rho_0-\omega\|_{C^{1, \alpha}}.
\]
Moreover, for any $t\in [t_0/2, t_0]$, and $k\geq 2$,
\[
\|\rho-\omega\|_{C^{k, \alpha}}\leq C_k \|\rho_0-\omega\|_{C^{1, \alpha}}
\]
\end{prop}
We assume $\delta_0\in (0, 1)$ is sufficiently small and denote, for $k\geq 2$,  \[V_0=\{\rho\in [\omega]: \|\rho-\omega\|_{C^{1, \alpha}}<\delta_0, \|\rho-\omega\|_{C^{k, \alpha}}\leq C_k\}.\] In particular for any $\rho\in V_0$, its potential volume $u$ is closed to $1$, such that the flow \eqref{conformal-flow1} is uniformly parabolic. Now we fix the constants $\delta_0, t_0, C_0, C_k$ and $V_0$. For any $t\in [t_0/2, t_0]$, $\rho(t)\in V_0$.  Since $\rho(t_0)\in V_0$,  we can get a solution of \eqref{conformal-flow1} over $[t_0, 2t_0]$ by applying Proposition \ref{short2} with initial date $\rho(t_0)$. For $t\in [t_0, 2t_0]$, 
we have 
\[
\|\rho-\omega\|_{C^{1, \alpha}}\leq C_0 \|\rho(t_0)-\omega\|_{C^{1, \alpha}} \leq C_0^2 \|\rho_0-\omega\|_{C^{1, \alpha}}
\]
and also 
\[
\|\rho-\omega\|_{C^{k, \alpha}}\leq C_k\|\rho(t_0)-\omega\|_{C^{k, \alpha}}\leq C_k^2 \|\rho_0-\omega\|_{C^{1, \alpha}}
\]
Now by choosing $\epsilon<<\delta_0$ sufficiently small, we can get that $\rho(t)\in V_0$ for $t\in [t_0, 2t_0]$ if $\|\rho_0-\omega\|_{C^{1, \alpha}}<\epsilon$. Applying this repeatedly, we can get finite time stability near $\omega$,
\begin{prop}[Finite time stability] For any $T\geq t_0$, there exists $\epsilon=\epsilon(T, \delta_0)$ such that if $\|\rho_0-\omega\|_{C^{1, \alpha}}<\epsilon$, the flow \eqref{conformal-flow1} exists over 
$[0, T]$ and $\rho(t)\in V_0$ for $t\in [t_0, T]$.  
\end{prop}

Now we are ready to proof Theorem \ref{stability}.
\begin{proof}
Fix $T_0>>1$ and the flow exists over $[0, T_0]$ for the initial data $\tilde \rho_0$ satisfying $\|\tilde\rho_0-\omega\|_{C^{1, \alpha}}<\epsilon_0=\epsilon(T_0, \delta_0)$ and $\rho(t)$ stays in $V_0$ over $[t_0, T_0]$. 
Now choose $\epsilon_1<<\epsilon_0$ such that the flow exists over $[0, T_1]$ with $T_1>>T_0$ and $\rho\in V_0$ over $[t_0, T_1]$, if the initial data satisfies
\[
\|\tilde \rho_0-\omega\|_{C^{1, \alpha}}<\epsilon_1. 
\]
 To obtain a uniform stability, we apply Proposition \ref{energydecay} to get, 
\[
\|\rho(t)-\omega\|^2_{L^2}\leq e^{-\frac{t}{c_2}} \|\rho_0-\omega\|^2_{L^2}\leq e^{-\frac{t}{c_2}}. 
\]
Hence we get that $\|\rho(t)-\omega\|_{C^{k, \alpha}}\leq C_k$ (since $\rho\in V_0$) and $\|\rho(t)-\omega\|^2_{L^2}\leq e^{-\frac{t}{c_2}}$, where $c_2$ is the Sobolev-Poincare constant. 
By the interpolation inequality, we get a uniform estimate, for $T_1$ is sufficiently large and $t\geq T_1/2$,
\begin{equation}\label{uniform}
\|\rho(t)-\omega\|_{C^{1, \alpha}}<\epsilon_0.
\end{equation}
In particular $\rho(T_1)$ satisfies 
\[
\|\rho(T_1)-\omega\|_{C^{1, \alpha}}<\epsilon_0.
\]
Now we fix $T_1$ and take $\epsilon=\epsilon_1$. 
Extend the flow with initial data $\rho(T_1)$, then we obtain the flow solution $\rho(t)$ over $[0, T_1+T_0]$ with $\rho(t)\in V_0$ over $[t_0, T_1+T_0]$. 
Since the uniform estimate \eqref{uniform} on $\rho(t)$ holds for all $t\geq T_1/2$, by repeating the argument
the flow exists over $[0, \infty)$ with uniformly bounded $C^{k, \alpha}$ norm and  the energy decay implies that $\rho(t)$ converges to $\omega$ in the exponential rate. 
\end{proof}

\subsection{Shi-type estimates}It is not clear whether the flow will exist for all time or it will develop a finite time singularity in general. To understand behavior of a possible finite time singularity, we develop the Shi-type estimate \cite{Shi} for the flow. 
Firstly we compute the evolution equations.
The evolution equation of $\rho$ reads, in local coordinates,
\[
\p_t \rho_{ij}=\frac{\Delta \rho_{ij}+(\cR(\rho))_{ij}}{\sqrt{u}}+\frac{u_i \rho_{jl, l}-u_j \rho_{il, l}}{2\sqrt{u^3}}
\]
In this section we assume that \begin{equation}\label{shi01}
\frac{1}{\sqrt{u}}+|\nabla \rho|^2\leq C_0.\end{equation}
Since along the flow $\|\rho\|_{L^2(M)}$ is uniformly bounded, the assumption \eqref{shi01} implies that $|\rho|$ is uniformly bounded. 
To derive the Shit-type estimates, we only need the following structures
\[
\p_t \rho=\frac{\Delta \rho+\cR\circ \rho}{\sqrt{u}}+\frac{\rho\circ \nabla \rho \circ\nabla \rho}{\sqrt{u^3}},
\]
where $A\circ B$ describe possible contraction of $A$ and $B$. 
Note that we understand \[\nabla u=\frac{1}{2} \langle (*\rho),  \nabla \rho\rangle =\rho\circ \nabla \rho.\]
We have the following,
\begin{prop}
Assuming \eqref{shi01}, we have
\[
\p_t |\nabla \rho|^2\leq \frac{\Delta |\nabla \rho|^2-\frac{3}{2}|\nabla^2 \rho|^2}{\sqrt{u}}+C,
\]
where $C$ is a uniformly bounded constant, depending on the metric $g$, its curvature and derivative of curvature, and the bound in \eqref{shi01}. More generally we have, for $k\geq 2$, 
\begin{equation}\label{shi02}
\p_t |\nabla^k \rho|^2\leq \frac{\Delta |\nabla^k \rho|^2-\frac{3}{2}|\nabla^{k+1} \rho|^2}{\sqrt{u}}+C_k \left(\sum_{i=2}^k |\nabla^i\rho|^{2\frac{(k+1)}{i}}+1\right).
\end{equation}
\end{prop}
\begin{proof}
This is a standard practice. With the symbolic notations, we have
\[
\p_t \nabla \rho=\frac{\Delta \nabla \rho}{\sqrt{u}}+\frac{\cR\circ \nabla \rho+\nabla \cR\circ \rho}{\sqrt{u}}+\frac{\nabla^2 \rho\circ \nabla \rho\circ \rho+\cR\circ \rho\circ\rho\circ\nabla \rho}{\sqrt{u^3}}+\frac{\rho\circ\rho\circ\nabla \rho\circ\nabla \rho\circ\nabla \rho}{\sqrt{u^5}}
\]
It follows that
\[
\p_t |\nabla \rho|^2\leq\frac{\Delta |\nabla\rho|^2-2|\nabla^2\rho|^2}{\sqrt{u}}+C(|\nabla^2\rho|+1)\leq \frac{\Delta |\nabla\rho|^2-\frac{3}{2}|\nabla^2\rho|^2}{\sqrt{u}}+C. 
\]
Note that with the assumption $|\nabla u|$ is bounded. Hence we have
\[
\p_t \nabla^2 \rho\leq \frac{\Delta \nabla^2 \rho}{\sqrt{u}}+C (|\nabla^3\rho|+|\nabla^2\rho|^2+|\nabla^2 \rho|+1)
\]
And hence 
\[
\begin{split}
\p_t |\nabla^2\rho|^2\leq& \frac{\Delta |\nabla^2 \rho|^2-2|\nabla^3 \rho|^2}{\sqrt{u}}+C|\nabla^2\rho|(|\nabla^3\rho|+|\nabla^2\rho|^2+|\nabla^2 \rho|+1))\\
\leq &\frac{\Delta |\nabla^2 \rho|^2-\frac{3}{2}|\nabla^3 \rho|^2}{\sqrt{u}}+C(|\nabla^2\rho|^3+1)
\end{split}
\]
This proves \eqref{shi02} for $k=2$. For general $k$, an inductive computation gives that
\[
\p_t |\nabla^k\rho|^2\leq \frac{\Delta |\nabla^k \rho|^2-2|\nabla^{k+1} \rho|^2}{\sqrt{u}}+C|\nabla^k\rho|\left(\sum |\nabla^{i_1}\rho| |\nabla^{i_2}\rho|\cdots |\nabla^{i_p}\rho|+1\right)
\]
where the sum is taken such that $i_1+i_2+\cdots+i_p\leq k+2$, with $1\leq i_a\leq k+1$. By applying Young's inequality
\[
a_1 a_2\cdots a_k\leq \sum \frac{a_i^{p_i}}{p_i}, \sum p_i^{-1}=1
\]
we can get
\[
\p_t |\nabla^k\rho|^2\leq \frac{\Delta |\nabla^k \rho|^2-\frac{3}{2}|\nabla^{k+1} \rho|^2}{\sqrt{u}}+C_k \left(\sum_{i=2}^k |\nabla^i\rho|^{\frac{2(k+1)}{i}}+1\right)
\]
\end{proof}

Now we state the Shi-type estimate,

\begin{prop}Suppose \eqref{shi01} holds over $[0, T]$, then we have the uniform bound for $\nabla^k\rho$ for  $t\in (0, T]$, $k\geq 2$, 
\begin{equation}\label{shi03}
|\nabla^k\rho|\leq C_k t^{-\frac{(k-1)}{2}}
\end{equation}
Such an estimate can be made local. Suppose \eqref{shi01} holds on $B_r(p)\times [0, T]$, where $B_r(p)$ is the geodesic ball of the background metric $g$, then over 
$B_{r/2}(p)\times [T/2, T]$, we have
\begin{equation}\label{shi04}
|\nabla^{k}\rho|\leq C_k(g, C_0, k, T, r). 
\end{equation}
\end{prop}

\begin{proof}Given \eqref{shi02}, the proof of \eqref{shi03} follows almost line by line  as in He-Li \cite[Proposition 2.1]{HeLi}.
The local estimate is also standard. 
 Denote $Q=(\mu+|\nabla \rho|^2)|\nabla^2\rho|^2$ with $\mu$ determined later.
We compute
\[
\begin{split}
\p_t Q=&\p_t |\nabla\rho|^2 \;|\nabla^2\rho|^2+(\mu+|\nabla\rho|^2) \p_t |\nabla^2 \rho|^2\\
\leq & \left(\frac{\Delta |\nabla \rho|^2-\frac{3}{2}|\nabla^2 \rho|^2}{\sqrt{u}}+C\right) |\nabla^2\rho|^2\\&+(\mu+|\nabla\rho|^2)\left(\frac{\Delta |\nabla^2 \rho|^2-\frac{3}{2}|\nabla^3 \rho|^2}{\sqrt{u}}+C(|\nabla^2\rho|^3+1)\right)\\
\leq& \frac{\Delta Q}{\sqrt{u}}-\frac{3}{2\sqrt{u}} |\nabla^2\rho|^4-\frac{3\mu}{2\sqrt{u}} |\nabla^3\rho|^2 +C_1(|\nabla^3\rho|\;|\nabla^2\rho|^2+|\nabla^2\rho|^3+1)
\end{split}
\]
Since $u\leq C_2=C_2(C_0, g)$ is uniformly bounded by our assumption.
By choosing $\mu$ sufficiently large, depending on $C_1$ and $C_2$, we have
\begin{equation}\label{shi05}
\p_t Q-\frac{\Delta Q}{\sqrt{u}}\leq -A_0 Q^2+A_1
\end{equation}
With \eqref{shi05} the localization is standard as well. Choose a cut-off function $\phi$ supported in $B_r(p)$, and $\phi=1$ in $B_{r/2}$. Denote, for a suitable constant $A_2=A_2 (r^{-1}, A_0, A_1, C_0)$, 
\[
H=\frac{A_2}{\phi^2}+\frac{1}{A_0 t}+\sqrt{A_1 A_0^{-1}}. 
\]
We compute
\[
\p_t H=-\frac{1}{A_0 t^2}, \Delta H=\frac{A_2(-2\phi \Delta \phi+6|\nabla \phi|^2)}{\phi^4}, H^2>\frac{A_2^2}{\phi^4}+\frac{1}{A_0^2 t^2}+\frac{A_1}{A_0}
\]
For $x\in B_r(p)$, we have $|\Delta \phi|+|\nabla \phi|^2\leq 10 r^{-2}$. Hence for $A_2$ sufficiently large, 
\[
\p_t H-\frac{\Delta H}{\sqrt{u}}>-A_0 H^2+A_1. 
\]
A maximum principle argument then implies that $Q\leq H$ over $B_r(p)\times [0, T]$. This proves \eqref{shi04} for $k=2$. For general $k$ take $Q_k=(\mu_k+|\nabla^{k-1}\rho|^2)|\nabla^k \rho|^2$ for sufficiently large $\mu_k$. 
\end{proof}

Hence by the Shi-type estimate, if the flow exists for the maximal interval $[0, T)$ with $T<\infty$, then it develops finite singularities at $t=T$ with
$\frac{1}{\sqrt{u}}+|\nabla \rho|^2$ blowing up. Suppose  $\frac{1}{\sqrt{u}}$ is bounded, only the Sobolev bound $\rho \in W^{1, p}$ for $p>4$ is needed ($\rho$ would be Holder continuous in particular). This is a well-known technique in the study of harmonic map theory. Even though the argument is by far standard for specialists, it is technically very involved. Since we do not need this result, we will not pursue such a generalization. 

\subsection{Long time existence and convergence in $W^{1, 2}$ assuming the bound of $|\nabla \log u|$}
In this section we assume $|\nabla \log u|\leq C_0$. Note that we assume the normalized condition
\[
\int_M u\text{dvol}=\int_M \text{dvol}=1.
\]
Hence  for any $t$ there exists $x_t\in M$ such that $u(x_t, t)=1$. If $|\nabla \log u|\leq C_0$, then there exists $C_1=C_1(C_0, g)$ such that
\[
-C_1\leq \log u\leq C_1
\]
This in particular implies that
\[
e^{-C_1}\leq u\leq e^{C_1},\; \text{and}\; |\nabla u|\leq C_0 e^{C_1}.
\]
We have the following,
\begin{prop}\label{extension}
Suppose the flow exists on $[0, T]$ with $T<\infty$ and $|\nabla \log u|\leq C_0$ over $[0, T]$. Then the flow can be extended across $T$. \end{prop}

\begin{proof}
We need to be a bit more precise here given the extra assumptions on $u$ and $\nabla u$.
We write the equation as
\begin{equation}
\p_t \rho=\frac{\Delta \rho}{\sqrt{u}}+\frac{\cR\circ \rho}{\sqrt{u}}+\nabla (u^{-\frac{1}{2}})\circ \nabla \rho
\end{equation}
We compute
\begin{equation}
\p_t |\rho|^2=\frac{\Delta |\rho|^2-2|\nabla \rho|^2}{\sqrt{u}}+\frac{\langle \cR (\rho), \rho\rangle}{\sqrt{u}}+\frac{\nabla \log u \circ \nabla \rho \circ \rho}{\sqrt{u}}
\end{equation}
Suppose $|\nabla \log u|\leq C_0$, then we have
\begin{equation}
\p_t |\rho|^2-\frac{\Delta |\rho|^2}{\sqrt{u}}\leq \frac{-2|\nabla \rho|^2}{\sqrt{u}} +C(|\rho|^2+|\nabla \rho| |\rho|)
\end{equation}
It follows that  $|\rho|^2\leq |\rho_0|^2 e^{Ct}$ at time $t$. 
A straightforward computation gives
\begin{equation}
\p_t |\nabla \rho|^2\leq \frac{\Delta |\nabla \rho|^2-2 |\nabla^2\rho|^2}{\sqrt{u}}+C(|\nabla^2 \rho| |\nabla \rho|+|\nabla \rho|^2+|\nabla \rho| |\rho|+|\nabla^2 u||\nabla \rho|)
\end{equation}
For the potential $u$, we compute
\begin{equation}
\p_t u=\frac{\Delta u}{\sqrt{u}}-\frac{1}{\sqrt{u}} \langle \nabla (*\rho), \nabla \rho \rangle+\frac{\nabla u\circ \nabla \rho \circ \rho}{\sqrt{u^3}}+\frac{\langle \cR(\rho), *\rho\rangle}{\sqrt{u}}
\end{equation}
A straightforward computation also gives, with uniform bounds on $u$ and $|\nabla u|$, 
\begin{equation}
\begin{split}
\p_t |\nabla u|^2&\leq  \frac{\Delta |\nabla u|^2-2|\nabla^2 u|^2}{\sqrt{u}}\\
&+C(|\nabla^2 u|+|\nabla^2\rho||\nabla \rho|+|\nabla^2 u||\nabla \rho|+|\nabla^2\rho|+|\nabla \rho|^2+|\nabla \rho|+|\rho|^2)
\end{split}
\end{equation}
Denote $f=|\nabla \rho|^2+A |\nabla u|^2+B|\rho|^2+1$ with $1<<A<<B)$, where the constants $A, B$ depend on the bound on $u$ and $\nabla u$. Then we have, for a uniformly bounded constant $C_1$, 
\begin{equation}\label{structure-1}
\p_t f-\frac{\Delta f}{\sqrt{u}}\leq  C_1 f.
\end{equation}
A straightforward maximum principle argument implies that 
\[
|\nabla \rho|^2\leq C_2 e^{C_1 t}
\]
Hence if the flow exists over $[0, T]$ for $T<\infty$ such that $|\nabla \log u|\leq C_0$. Then $|\nabla \rho|^2\leq C_2 e^{C_1 T}$. By the Shi-type estimate, we know that the flow can be extended across $T$. In other words, the flow exists as long as $|\nabla \log u|$ stays bounded. 
\end{proof}

We have the following convergence result when assuming $|\nabla \log u|\leq C_0$. 
\begin{prop}\label{convergence01}Suppose the flow exists for all time with $|\nabla \log u|\leq C_0$, for an initial symplectic form $\rho_0\in [\omega_0]$. Then the flow converges in $L^2$ to $\omega_0$ when $t\rightarrow \infty$. 
\end{prop}

\begin{proof}If one obtains a uniform estimate of $|\nabla \log u|\leq C_0$ along the flow, then the long time existence follows Proposition \ref{extension}. 
The point bounds on $|\rho|^2, |\nabla \rho|^2$ etc. derived in the last section assuming $|\nabla \log u|\leq C_0$ depend on a time factor $e^{Ct}$, which is not bounded when $t\rightarrow \infty$. 
To derive convergence, we use the $L^2$ estimate. 
Recall that we have
\begin{equation}
\p_t \int_M |\rho-\omega_0|^2 =-\int_M \frac{|d^*\rho|^2}{\sqrt{u}}\leq -C_1 \int_M |d^*\rho|^2
\end{equation}
We compute, using $|\nabla \log u|\leq C_0$, 
\begin{equation}
\begin{split}
\p_t \int_M |d^*\rho|^2=&-\int_M \left\langle d^*\rho, d^*  d\left( \frac{d^*\rho}{\sqrt{u}}\right)\right\rangle=-\int_M \left\langle dd^*\rho,  d\left( \frac{d^*\rho}{\sqrt{u}}\right)\right\rangle\\
\leq &-\int_M \frac{|dd^*\rho|^2}{\sqrt{u}}+ C\int_M |dd^*\rho| |d^*\rho|.
\end{split}
\end{equation}
Denote 
\[
Q_1=\int_M |d^*\rho|^2+A_1  \int_M |\rho-\omega_0|^2.
\]
It follows that, for $A_1$ sufficiently large, 
\begin{equation}\label{key-01}
\p_t Q_1\leq -C_1 \left(\int_M |dd^*\rho|^2+\int_M |d^*\rho|^2\right)
\end{equation}
An application of the Poincare inequality implies that $\p_t Q_1\leq - C Q_1$ for some positive constant $C$. 
It follows that for $t\in [0, \infty)$
\begin{equation}\label{abc1}
\int_M |d^*\rho|^2+A_1  \int_M |\rho-\omega_0|^2\leq A_0 e^{-Ct}.
\end{equation}
When $t\rightarrow \infty$, $\rho(t)$ has uniformly bounded $W^{1, 2}$ norm, and hence converges strongly in $L^2$ and weakly in $W^{1, 2}$ to a limit form $\rho_\infty$ with $\rho_\infty\in [\omega_0]$,  $d\rho_\infty=0$, $\rho_\infty\in W^{1, 2}$.
The convergence implies that
\[
\int_M |d^*\rho_\infty|^2\leq \liminf \int_M |d^*\rho|^2=0, \int_M |\rho_\infty-\omega_0|^2=0.
\]
It follows that $\rho_\infty=\omega_0$. 
\end{proof}

\subsection{Smooth convergence}
To obtain smooth convergence, we need to obtain the pointwise bound on $|\rho|$ and $|\nabla \rho|$. This requires more careful analysis. We have the following, 
\begin{prop}Suppose $|\nabla \log  u|\leq C_0$ over $[0, \infty)$ along the flow. Then $|\rho|$ and $|\nabla \rho|$ are both uniformly bounded. It follows that the flow converges to $\omega$ smoothly when $t\rightarrow \infty$ with uniform control. 
\end{prop}

\begin{proof}The main point is to bound $f=|\nabla \rho|^2+A|\nabla u|^2+B|\rho|^2+1$, which satisfies \eqref{structure-1}, given the uniform control on $|\nabla \log u|$. 
We recall the inequality
\begin{equation}
\p_t f\leq \frac{\Delta f}{\sqrt{u}}+C_1 f. 
\end{equation}
We apply the parabolic Moser iteration, which is a standard and very powerful tool for specialists. Since the argument is involved  and technically delicate in its nature, we include the details for completeness. For simplicity we work on the time interval $[0, 1]$. 
For $2p-1>0$ (say take $p\geq \frac{5}{8}$),  we compute
\begin{equation}
\begin{split}
\p_t \int_M f^{2p}\leq& 2p\int_M f^{2p-1}\frac{\Delta f}{\sqrt{u}}+2pC_1\int_M f^{2p}\\
=&-2p(2p-1)\int_M f^{2p-2} \frac{|\nabla f|^2}{\sqrt{u}}-2p\int_Mf^{2p-1}\nabla f \nabla  (u^{-\frac{1}{2}})+2p\int_M f^{2p}\\
\leq&-4c_0 p(2p-1)\int_M f^{2p-2}|\nabla f|^2+pC\int_M f^{2p-1} |\nabla f|+2p\int_M f^{2p}\\
\leq &-4c_0\left(2-\frac{1}{p}\right)\int_M |\nabla (f^p)|^2+C\int_M f^{p-1}|\nabla (f^p)|+2p\int_M f^{2p}\\
\leq &-c_0\int_M |\nabla (f^p)|^2+pC \int_M f^{2p},
\end{split}
\end{equation}
where we use the Holder inequality and Young's inequality. We use  $c_0$ to denote a uniform positive constant depending on $\min u^{-\frac{1}{2}}$, and use $C$ to denote a uniformly bounded positive constant which can vary line by line. Applying a Sobolev inequality to $h=f^p$
\[
\left(\int_M h^4\right)^{\frac{1}{2}}\leq C_0 \int_M (|\nabla h|^2+h^2)
\]
we have the following inequality, for $p\in \left[\frac{5}{8}, \infty\right)$, 
\begin{equation}
\p_t \int_M f^{2p}+c_1\left(\int_M f^{4p}\right)^{\frac{1}{2}}\leq pC \int_M f^{2p}
\end{equation}
Now given $0\leq a<b\leq 1$,  let $\phi: [0, 1]\rightarrow [0,1]$ be a Lipschitz function such that \[\phi(s)=0, s\in [0, a];\, \phi(s)=1, s\in [b, 1]\] and $\phi(s)$ is linear for $s\in [a, b]$.  Then  \[\phi^{'}(s)=\frac{1}{b-a} ,s\in (a, b);\, \phi^{'}=0, s\in [0, a)\cup (b, 1].\]
Taking integration with respect to $t$ over $[0, t_0]$ for any $t_0\in [a, 1]$, we have
\begin{equation}
\int_0^{t_0} \phi^2(t)\p_t\int_M f^{2p}+c_0 \int_0^{t_0} \phi^2(t) \left(\int_M f^{4p}\right)^{\frac{1}{2}}\leq pC\int_0^{t_0} \phi^2(t)\int_M f^{2p}.
\end{equation}
Integration by parts on the first term above, we have
\begin{equation}\label{structure-2}
\left(\phi^2\int_M f^{2p}\right) (t_0)+c_0\int_0^{t_0} \phi^2 \left(\int_M f^{4p}\right)^{\frac{1}{2}}\leq \left(pC+\frac{2}{b-a}\right)\int_a^{t_0}\int_M f^{2p}\end{equation}
 We have, by \eqref{structure-2},
\begin{equation}\label{structure-3}
\max_{[b, 1]}\int_M f^{2p}\leq  \left(pC+\frac{2}{b-a}\right)\int_a^1\int_M f^{2p}
\end{equation}
and 
\begin{equation}
\int_{b}^1\left(\int_M f^{4p}\right)^{\frac{1}{2}}\leq \left(pC+\frac{2}{b-a}\right)\int_a^1\int_M f^{2p}
\end{equation}\label{structure-4}
 Now for any $t\in [b, 1]$, applying the H\"older inequality we get
 \begin{equation}
 \int_M f^{3p} \leq \left(\int_M f^{4p}\right)^{\frac{1}{2}}\left(\int_M f^{2p}\right)^{\frac{1}{2}}\leq \left(\max_{[b, 1]}\int_M f^{2p}\right)^{\frac{1}{2}}\left(\int_M f^{4p}\right)^{\frac{1}{2}}
 \end{equation}
 It follows that over the space time $[b, 1]\times M$, 
 \begin{equation}\label{structure-5}
 \int_b^1\int_M f^{3p}\leq \left(\max_{[b, 1]}\int_M f^{2p}\right)^{\frac{1}{2}} \int_b^1 \left(\int_M f^{4p}\right)^{\frac{1}{2}}.
 \end{equation}
 Applying \eqref{structure-3} and \eqref{structure-4} to \eqref{structure-5}, we have
 \begin{equation}
\left( \int_b^1\int_M f^{3p}\right)^2\leq \left(pC+\frac{2}{b-a}\right)^3\left(\int_a^1\int_M f^{2p}\right)^3
 \end{equation}
 We write the above as the following,
 \begin{equation}
 \|f\|_{L^{3p}([b, 1]\times M)}\leq \left(pC+\frac{2}{b-a}\right)^{\frac{1}{2p}} \|f\|_{L^{2p}([a, 1]\times M)}. 
 \end{equation}
 In other words, we have,
 \begin{equation}\label{structure-6}
 \log  \|f\|_{L^{3p}([b, 1]\times M)}\leq \frac{1}{2p} \log \left(pC+(b-a)^{-1}\right)+\log \|f\|_{L^{2p}([a, 1]\times M)}. 
 \end{equation}
 Now take a sequence of $p_i$ and $b_i$ such that
 \[
 p_i=\frac{1}{2}\left(\frac{3}{2}\right)^{i+1}, b_i=\frac{1}{2}-\frac{1}{2^{i+1}}, i=0, 1, 2, \cdots
 \]
 Denote $q_i=2p_{i}$ as well.
 For each $i$, applying \eqref{structure-6} with $p=p_i$, $a=b_i, b=b_{i+1}$, we get
 \begin{equation}\label{structure-7}
 \log \|f\|_{L^{q_{i+1}}([b_{i+1}, 1]\times M)}\leq \left(\frac{2}{3}\right)^{i+1}\log (2^{i+2}+Cq_i)+\log \|f\|_{L^{q_i}([b_i, 1]\times M}).
 \end{equation}
 By letting $i\rightarrow \infty$ and the fact that
 $\sum_i \left(\frac{2}{3}\right)^i i<\infty$,
\eqref{structure-7} gives 
 \begin{equation}\label{key-02}
 \log \|f\|_{L^\infty([\frac{1}{2}, 1]\times M)}\leq C+\log \|f\|_{L^{\frac{3}{2}}([0, 1]\times M)}
 \end{equation}
 By choosing the constants in a more careful way, we can indeed get that, for any $r>1$, since \eqref{structure-6} holds only for $2p>1$, 
\begin{equation}
 \log \|f\|_{L^\infty([\frac{1}{2}, 1]\times M)}\leq C_1(r)+\log \|f\|_{L^{r}([0, 1]\times M)},
 \end{equation}
where the constant $C_1$ depends on $r$, and it could blow up when $r\rightarrow 1$. To overcome this difficulty, we prove that $f$ is indeed in $L^{\frac{3}{2}}([t, t+1]\times M)$ for any $t\geq 0$. Since $\rho$ is in $W^{1, 2}(M)$ and hence $|\rho|$ is in $L^4(M)$. Note that, given $|\nabla u|$ bounded, 
\[
\int_M f^{\frac{3}{2}}\leq 2\int_M |\nabla \rho|^3+C \int_M |\rho|^3+C.
\]
Hence to show that $f\in L^{\frac{3}{2}}([t, t+1]\times M)$ for each $t$, it suffices to show that for any $t$
\begin{equation}\label{structure-8}
\int_t^{t+1}\int_M |\nabla \rho|^3\leq C. 
\end{equation}
At the moment we have no uniform control on 
$\int_M |\nabla \rho|^3$ for each time slice.  To obtain \eqref{structure-8}, we need to use \eqref{key-01}. Taking integral over $[t, t+1]$, we get that
\[
\int_t^{t+1} \int_M |d^*d\rho|^2+\int_t^{t+1} \int_M |d^*\rho|^2\leq Q_1(t)-Q_1(t+1)\leq C.
\]
Using the  Bochner-Weitzenbock formula we obtain that
\[
\int_t^{t+1} \int_M |\nabla^2\rho|^2+\int_t^{t+1} \int_M |\nabla\rho|^2\leq C.
\]
By applying Sobolev inequality on each time slice, we obtain 
\[
\left(\int_M |\nabla \rho|^4\right)^{\frac{1}{2}}\leq  C\left(\int_M |\nabla^2\rho|^2+|\nabla \rho|^2\right)
\]
Hence we have 
\[
\int_t^{t+1} \left(\int_M |\nabla \rho|^4\right)^{\frac{1}{2}}\leq C.
\]
Now we have, since $|\nabla \rho|\in L^2(M)$, 
\[
\int_t^{t+1}\int_M |\nabla \rho|^3\leq \int_t^{t+1} \left(\int_M |\nabla \rho|^4\right)^{\frac{1}{2}} \left(\int_M |\nabla \rho|^2\right)^{\frac{1}{2}}\leq C \int_t^{t+1}\left(\int_M |\nabla \rho|^4\right)^{\frac{1}{2}}\leq C. 
\]
This prove \eqref{structure-8}. Applying \eqref{key-02} over $[t, t+1]$ and \eqref{structure-8},  we have
the uniform control on $f$ over $[t+\frac{1}{2}, t]$ for any $t$. This give the uniform bound on $|\rho|$ and $|\nabla \rho|$, which implies the smooth convergence of $\rho\rightarrow \omega_0$ along the flow. 
\end{proof}

The assumption $|\nabla \log u|\leq C$ can be weaken considerately, by assuming a uniform control $u+u^{-1}+\|\nabla u\|_{L^p}$ for some $p$ sufficiently large.
In particular, we can ask the following technical question,

\begin{pl}
For long time existence, is it sufficient to assume a uniform control on $u+u^{-1}+\|\nabla u\|_{L^p}$, for $p>4$? for $p=4$?
\end{pl}
The problem above seems plausible, by refining our argument carefully. But it would need substantial technical improvement to assume  $u+u^{-1}+\|\nabla u\|_{L^4}\leq C$. We will leave these technical issues for further investigation. A more ambitious problem would be only assuming $L^\infty$ control of $|\rho|$ and nondegeneracy of $u$. 
\begin{pl}Does the flow extend by only assuming a uniform control of $u^{-1}+|\rho|^2$?
\end{pl}

\subsection{Examples and discussions}\label{example}
In this section we discuss some examples. 
Firstly we consider an example on $T^4=S^1\times S^1\times S^1\times S^1$. Let $\omega_0= dx_0\wedge dx_1+dx_2\wedge dx_3$ and $g_0=\sum dx_i\otimes dx_i$. Suppose $\rho_0=u_0(x_1, x_2) dx_1\wedge dx_2+dx_3\wedge dx_4$ for a positive function $u_0: T^2=S^1\times S^1\rightarrow \R_{+}$ satisfying 
\begin{equation}\label{torus1}
\int_{T^2} u_0 dx_1\wedge dx_2=\int_{T^2} dx_1\wedge dx_2. 
\end{equation}
Clearly $d\rho_0=0$ and \eqref{torus1} implies that $\rho_0\in [\omega_0]$. A straightforward computation implies that the flow \eqref{conformal-flow1} reduces to the parabolic flow of $u(t): T^2\rightarrow \R_{+}$, 
\[
\p_t u=\frac{\Delta u}{\sqrt{u}}-\frac{|\nabla u|^2}{2\sqrt{u^3}}
\]
Note that one can write the equation as
\begin{equation}\label{torus2}
\p_t u=2\Delta \left({u^{\frac{1}{2}}}\right)
\end{equation}
This is an example of \emph{fast diffusion equation}. Since $u_0$ is assumed to be positive, the maximum principle implies that $u$ stays uniformly bounded and uniformly positive. By the Krylov-Safanov estimate, $u$ has uniformly bounded H\"older norm. Hence \eqref{torus2} has a unique long time smooth solution with uniformly bounded H\"older norms by the standard parabolic theory. The gradient flow nature implies that $u$ converges smoothly to $u_\infty\equiv 1$. 

This extends to the example $(M^4, \omega_0=\omega_1\oplus \omega_2)=(\Sigma_1, \omega_1)\times (\Sigma_2, \omega_2)$ for two closed surfaces with two volume forms $\omega_1$, $\omega_2$.
Suppose $\rho_0\in [\omega_0]$ of the form $\rho_0=u_0\omega_1\oplus \omega_2$ for $u_0: \Sigma_1\rightarrow \R_{+}$ such that 
\[
\int_{\Sigma_1} u_0 \omega_1=\int_{\Sigma_1} \omega_1. 
\]
The flow \eqref{conformal-flow1} is then reduced to the fast diffusion equation in \eqref{torus2} in the classical setting. The computation is almost identical. 
The example above does not extend to the case when $\rho_0=v_0(x_1, x_2) dx_1\wedge dx_2+w_0(x_3, x_4) dx_3\wedge dx_4$ satisfying $v_0, w_0>0$, and 
\[
\int_{T^2} v_0 dx_1\wedge dx_2=\int_{T^2} w_0 dx_3\wedge dx_4=\int_{T^2} dx_1\wedge dx_2=\int_{T^2} dx_3\wedge dx_4. 
\]
In this case $\rho_0\in [\omega_0]$, but the product structure is not preserved along the flow \eqref{conformal-flow1}. Note that in this example, $\rho_0$ is isotropic to $\omega_0$ in an obvious way. 

\begin{pl}
Does the flow \eqref{conformal-flow1} exist for all time and converge to $\omega_0$ on $T^4$ in the above example, with the initial symplectic form \[\rho_0=v_0(x_1, x_2) dx_1\wedge dx_2+w_0(x_3, x_4) dx_3\wedge dx_4?\]
\end{pl}

The product example above might give some indication of how to choose the metric $\tilde g=g(Q, \cdot)$ and the corresponding Hodge star $\tilde *$. Suppose the matrices $A$ and  $B$ are defined as
$\rho=g(A\cdot, \cdot)$ and $*\rho=g(B\cdot, \cdot)$. For $g_0$ and $\rho_0$ given above on $T^4$, we have 
\[
A=\begin{pmatrix}
0 & v & 0 &0\\
-v&0 &0&0\\
0&0 &0 &w\\
0&0&-w&0
\end{pmatrix},\,\,\, B=\begin{pmatrix}
0 & w & 0 &0\\
-w&0 &0&0\\
0&0 &0 &v\\
0&0&-v&0
\end{pmatrix}
\]
Hence $a=-A^2=\text{diag}(v^2, v^2, w^2, w^2)$ and $b=-B^2=\text{diag}(w^2, w^2, v^2, v^2)$. One can define $Q=f a^r$ or $fb^r$ for some positive function $f$ and some power $r$. For the conformal metric we have $Q=\sqrt{u} I$. For the metric $g_\rho$, $Q$ is given by $u^{-1}a$. Neither of the choice of $Q=\sqrt{u}I$ nor $Q=a$ (or $Q=u^{-1}a$) would preserve the product structure. In particular the flow becomes complicated for these choices. 
But if we choose $Q=b$ (or $b^r$ for $r>0),$ then the flow preserves the product structure such that
\[
\p_t v=\Delta \left(-\frac{1}{v}\right);\; \p_t w=\Delta \left(-\frac{1}{w}\right)
\]
If we choose $Q=\sqrt{b}$, then the flow equation preserves the product structure such that
\[
\p_t v=\Delta (\log v);\; \Delta w=\Delta (\log w)
\]
For such a choice, the flow equation becomes the classical fast diffusion equation and the flow exists for time and converges to $\omega_0$ on $T^4$ by the standard theory of parabolic equations. 

\begin{pl}Let $(M, \omega, g)$ be an almost K\"ahler manifold of dimension four. For each symplectic form $\rho\in [\omega]$.  We define the metric $\tilde g=g(b\cdot, \cdot)$.  Consider the flow 
\[
\p_t \rho=d(\tilde * d*\rho)
\]
Does the flow exist for all time?
\end{pl}
We will analyze the different Hodge flows in the next section. 
For the second example, we consider $\rho_0\in [\omega_0]$ on $(T^4, \omega_0, g_0)$, with 
\[
\rho_0=dx_1\wedge dx_2+dx_3\wedge dx_4+d\theta,\; \text{with}\; \theta=a(x_1, x_2)dx_3+b(x_1, x_2) dx_4
\]
where $a, b$ are two functions on $S^1\times S^1$ (hence  periodic functions on $\R\times \R$). The nondegenerate condition implies that
\[
u_0=1-a_{x_1} b_{x_2}+a_{x_2} b_{x_1}>0. 
\]
Note that for this example, $\rho_0$ is isotropic to $\omega_0$ via $\rho_s= \omega_0+s d\theta$ for $s\in [0, 1]$, with $\rho_s$ nondegenerate for each $s$. 
The linear Hodge Laplacian flow reduces to the heat equations for $a, b$ separately
\[
\p_t a=\Delta a=a_{x_1x_1}+a_{x_2x_2};\;\;\;
\p_t b=\Delta b=b_{x_1x_1}+b_{x_2x_2}
\]
However along the Hodge Laplacian flow, the positivity condition on $u=1-a_{x_1} b_{x_2}+a_{x_2} b_{x_1}$ might not hold. 
This fact should be well-known for experts. Since we cannot find a direct reference, we construct an explicit example as follows. 

\begin{exam}\label{negative-01}
We consider the linear Hodge flow for $\rho=dx\wedge dy+dz\wedge dw+d(a(x, y)dz+b(x, y)dw)$ on $T^4$. The flow is reduced to the heat equation on $T^2$
\[
\p_t a=\Delta a,\;\; \p_t b=\Delta b,
\]
where $a$ and $b$ are $2\pi$ periodic functions for $x$ and $y$ respectively. Take $a(x, y)=a(x)$ and suppose $a_x=f(x)$. While $b(x, y)=h(x) e(y)$. Along the linear flow, $f(x), h(x), e(y)$ are all periodic functions (functions on $S^1$) satisfying the heat equation
\[
\p_t f =f_{xx}, \p_t h=h_{xx},\;\text{and}\,\; \p_t\, e=e_{yy}. 
\]
Denote $e_y=c(y)$.  
For the initial data, $u_0(x, y)=1-f_0(x)h_0(x) c_0(y)$ and $u(x, y, t)=1-f(x, t) h(x, t) c(y, t)$, since $a(x, y, t)=f(x, t)$ is independent of $y$. 
There is a restriction on $f$ and $c$, that
\[
\int_0^{2\pi} f_0(x) dx=\int_0^{2\pi} c_0(y) dy=0
\]
We choose $f_0$ such that $f_0 = 0$ for $x\in [0, \pi]$ and $f_0(x)=\sin{2x}$, $x\in [\pi, 2\pi]$ and $h_0(x)=0$  for $x\in [\pi, 2\pi]$ and $h_0(x)=\sin(2x)$ for $x\in [0, \pi]$.  Hence $u_0(x, y)=1$ for any $(x, y)\in S^1\times S^1$. Let $f(x,t)$ and $h(x, t)$ solve the heat equation with initial datum $f_0$ and $h_0$ respectively. Clearly $f(x, t)h(x, t)$ is not identically zero for all $t>0$, and  $f(x, t)$ and $h(x, t)$ both converge to zero when $t\rightarrow \infty$. Let $\max f(x, 1) h(x, 1)= \frac{1}{A}>0$. Let $c(y, t)=A_0 e^{-t} \sin y$, then $u(x, y, 1)<0$ for some point if $A_0>A e$.  By choosing $A_0$ sufficiently large depending on $t_0$, $u(x, y, t_0)$ can become negative for arbitrarily small $t_0>0$. 
\end{exam}

The conformal flow \eqref{conformal-flow1} is reduced to the system
\begin{equation}\label{torus3}
\p_t a=\frac{\Delta a}{\sqrt{u}};\;\;\;\p_t b=\frac{\Delta b}{\sqrt{u}}
\end{equation}
and the factor $\frac{1}{\sqrt{u}}$ forces the positivity of $u$ before developing a singularity. 
\begin{prop}
Suppose $a(t), b(t): T^2\rightarrow \R$ satisfy \eqref{torus3}, then $\rho=\omega_0+d\theta$ with 
$\theta=a dx_3+bdx_4$
satisfies \eqref{conformal-flow1}.
\end{prop}

\begin{proof}This is a straightforward computation and the key point is that $d^*\theta=0$. 
Suppose $a, b$ are unique solutions of \eqref{torus3}. Then 
\[
\p_t \rho=d(\p_t \theta)=d\left(\frac{\Delta a\; dx_3+\Delta b\; dx_4}{\sqrt{u}}\right)
\]
Since $d^*\theta=0$, we have
\[-d*\rho=-d^* d\theta=\Delta \theta=\Delta a \;dx_3+\Delta b \; dx_4\]
Hence we have
\[
\p_t \rho=-d\left(\frac{d^*\rho}{\sqrt{u}}\right)
\]
\end{proof}

The system \eqref{torus3} is a quasilinear parabolic system for $(a, b)$ on $T^2$, together with a positivity restriction $u=1-a_{x_1} b_{x_2}+a_{x_2} b_{x_1}>0$. 
\begin{pl}Does \eqref{torus3} have a long time solution such that $\rho$ converges to $\omega_0$?
\end{pl}

 It is not clear at the moment whether the flow \eqref{conformal-flow1} would develop a singularity or not, when $n=4$. 
 The study of two special examples above might give some hint on this problem.

\subsection{Singularity and soliton}
When $n=6$, the corresponding geometric flow would have to develop singularities in general. In dimension $6$ we choose $r=1/3$ to consider the flow
\[
\p_t \rho=-d\left(\frac{d^*\rho}{u^{\frac{1}{3}}}\right)
\]
With almost identical arguments as in $n=4$,
one can show that if $|\nabla \log u|$ stays uniformly bounded along the flow, then the flow exists for all time and converges to $\omega$ strongly in $L^2$. (If one assume further that $|\nabla \rho|$ is bounded, then the convergence would be in $C^\infty$ topology).
One expects that in general the symplectic forms in a cohomology class $[\omega]$ is not path-connected. A particular example on $T^2\times S^2\times S^2$ is provided by McDuff. In this case, the flow has to develop singularities. The flow can be viewed as a gradient flow and the normalized energy decreases strictly along the flow (except at $\omega$), it would be interesting to analyze the formation of singularities along the flow. We pose a concrete question on formation of singularities, 
\begin{pl}Prove (or disprove) that along the conformal Hodge flow, there is no type-I singularity.
\end{pl}

One considers the solitons on $\R^4$. A particular type is as follows, 
\begin{equation}\label{soliton1}
d\left(\frac{d^*\rho_0}{\sqrt{u_0}}\right)=d(i_X\rho_0),
\end{equation}
where $X$ is a Killing vector field of the Euclidean metric. Let $\psi_t$ be the flow generated by $X$, then $\rho(t)=\psi_t^* \rho_0$ satisfies the flow equation on $\R^4$,
\[
\p_t \rho=-d\left(\frac{d^*\rho}{\sqrt{u}}\right)
\]
We can describe a family of special solutions to \eqref{soliton1} as follows. 
Let 
\[
\rho_0= a(x, y) dx\wedge dy+dw\wedge dz
\]
Choose $X=(V^1, V^2, V^3, V^4)$, a constant vector field. Then \eqref{soliton1} is reduced to
\begin{equation}\label{soliton2}
\left(\frac{a_x}{\sqrt{a}}\right)_x+\left(\frac{a_y}{\sqrt{a}}\right)_y+V^1 a_x+V^2 a_y=0.
\end{equation}                                    
It would be an interesting question to classify solutions to \eqref{soliton1}.

\section{Other nonlinear Hodge flows in dimension four}\label{example2}
In this section we discuss different choices of nonlinear Hodge flows. We write the flow equation as
\begin{equation}\label{hflow}
\p_t \rho_{ij}=(h_{ik}\rho_{kl, l})_{,j}-(h_{jk}\rho_{kl, l})_{,i}
\end{equation}
where $h_{ij}=Q_{ij} \sqrt{\det{Q}}^{-1}$ is a positive definite symmetric matrix such that $h_{ij}\rho_{jk}=\rho_{ik}h_{jk}$.  Note that $h_{ij}$ and $Q_{ij}$ are uniquely determined by each other. We fix the background structure $(M, \omega_0, g_0)$.  
For simplicity, we consider the case on  the flat four torus, which will simplify computations since the curvature terms all vanish. We write $g_0=\sum_i^4 d\theta_i^2$ and $\omega_0=d\theta_1\wedge d\theta_2+d\theta_3\wedge d\theta_4$. Suppose $\rho\in [\omega_0]$ is a symplectic form with 
\[
\rho=\sum_{i<j}\rho_{ij} d\theta_i\wedge d\theta_j=\sum_{i\neq j}\frac{1}{2}\rho_{ij}d\theta_i\wedge d\theta_j.
\]

\begin{prop}Along the flow, we have the following,
\begin{equation}
\begin{split}
\p_t |\rho|^2=&\rho_{ij}h_{ik}\Delta \rho_{kj}+2\rho_{ij}h_{ik,\,j}\rho_{kl,\,l}\\
\p_t u=&\frac{1}{2}(*\rho)_{ij}h_{ik}\Delta \rho_{kj}+(*\rho)_{ij}h_{ik, j}\rho_{kl, l}
\end{split}
\end{equation}
\end{prop}
\begin{proof}
We compute
\[
\begin{split}
\p_t |\rho|^2&=2\sum_{i<j} \rho_{ij}\p_t \rho_{ij}=\sum_{i\neq j} \rho_{ij} \p_t \rho_{ij}\\
&=2\rho_{ij}(h_{ik}\rho_{kl, l})_{j}\\
&=2\rho_{ij}h_{ik}\rho_{kl, lj}+2\rho_{ij} h_{ik, j} \rho_{kl, l}
\end{split}
\]
By $d\rho=0$, we have $\rho_{kl,\,lj}=\rho_{jl,\,kl}+\Delta {\rho_{kj}}$. Using $\rho_{ij}h_{ik}=\rho_{ki}h_{ij}$, we deduce \[2\rho_{ij}h_{ik}\rho_{kl,\,lj}=\rho_{ij}h_{ik}\Delta \rho_{kj}.\]
Similarly, we compute
\begin{equation*}
\p_t u=\frac{1}{2}(*\rho)_{ij}h_{ik}\Delta \rho_{kj}+(*\rho)_{ij}h_{ik, j}\rho_{kl, l}
\end{equation*}
This completes the computation. 
\end{proof}

First we consider the linear Hodge flow. 

\subsection{The linear Hodge flow $h_{ij}=\delta_{ij}$}
When $Q_{ij}=h_{ij}=\delta_{ij}$, the flow is the Hodge flow
\[
\p_t \rho=\Delta \rho
\]
The standard theory shows that the flow exists for all time with uniform control, and converges to $\omega_0$ at exponential rate. But the potential $u$ does not remain positive in general. In the last section we construct an example, for any given time $t_0>0$, one can choose an initial symplectic form $\rho_0$ with $u_0=1$, but along the flow $u$ can become negative at $t_0$. Nevertheless the linear Hodge flow can indicate some key difficulties in the study of nonlinear Hodge flow, in particular for conformal Hodge flows. 
We have
\[
\begin{split}
&\p_t |\rho|^2=\Delta |\rho|^2-2|\nabla \rho|^2, \p_t u=\Delta u-\langle\nabla \rho, \nabla (*\rho)\rangle\\
&\p_t |\rho^{+}|^2=\Delta |\rho^{+}|^2-2|\nabla \rho^{+}|^2, \p_t |\rho^{-}|^2=\Delta |\rho^{-}|^2-2|\nabla\rho^{-}|^2. 
\end{split}
\]
The eigenvalues of $A$ are $\pm \i \l_1, \pm \i \l_2$, which are determined by
\[
\l_1=\frac{1}{\sqrt{2}} (|\rho^{+}|+|\rho^{-}|), \l_2=\frac{1}{\sqrt{2}} (|\rho^{+}|-|\rho^{-}|)
\]
Note that $\l_1$ and $\l_2$ are smooth functions at points where $\rho^{-}\neq 0$. We compute
\[
\begin{split}
\p_t \l_1=\Delta \l_1-\frac{|\nabla \rho^{+}|^2-|\nabla |\rho^{+}||^2}{\sqrt{2}|\rho^{+}|}-\frac{|\nabla \rho^{-}|^2-|\nabla |\rho^{-}||^2}{\sqrt{2}|\rho^{-}|}\\
\p_t \l_2=\Delta \l_2-\frac{|\nabla \rho^{+}|^2-|\nabla |\rho^{+}||^2}{|\rho^{+}|}+\frac{|\nabla \rho^{-}|^2-|\nabla |\rho^{-}||^2}{|\rho^{-}|}
\end{split}
\]
We denote, at points where $\rho^{-}\neq 0$
\begin{equation}
J=\frac{|\nabla \rho^{+}|^2-|\nabla |\rho^{+}||^2}{\sqrt{2}|\rho^{+}|}, K=\frac{|\nabla \rho^{-}|^2-|\nabla |\rho^{-}||^2}{\sqrt{2}|\rho^{-}|}
\end{equation}
The maximum principle argument implies the larger eigenvalue $\l_1$ is bounded above, but fails to preserve the positivity of the smaller eigenvalue $\l_2$. (Strictly speaking, our computation only works at points where $\l_1\neq \l_2$, but with modifications the argument can work at points where $\l_1=\l_2$ as well). Note that at the minimum of $\l_2$, we have
\[
\nabla |\rho^{+}|=\nabla |\rho^{-}|.
\]
At the first time $\l_2=0$ and hence $|\rho^{+}|=|\rho^{-}|$ as well, the failure to preserve the positivity (nonnegativity) of $\l_2$ lies in the fact that at such points (the first time $\l_2=0$), the following control does not hold
\begin{equation}\label{u-1}
|\nabla \rho^{+}|\leq |\nabla \rho^{-}|
\end{equation}
Similar consideration of evolution equation on $u$ leads to the same observation, that at the point where $u=0$ for first time, \eqref{u-1} does not hold (it has to fail in general). At the maximum point of $u^{-1}|\rho|^2$, we have
\begin{equation}\label{u-1-2}
|\rho^{-}|^2 \nabla |\rho^+|^2=|\rho^{+}|^2\nabla |\rho^{-}|^2. 
\end{equation}
However $u^{-1} |\rho|^2$ fails to be bounded above along the flow, which shows that the following control at maximum point of $u^{-1}|\rho|^2$ fails as well (when $u$ is positive),
\begin{equation}\label{u-2}
|\rho^{-}|^2 |\nabla \rho^{+}|^2\leq |\rho^{+}|^2 |\nabla \rho^{-}|^2. 
\end{equation}
The failure of positivity of $u$ along the linear Hodge flow lies in the difference between \eqref{u-1-2} and \eqref{u-2} from technical perspective. 
Now we compute the evolution equations of $\l_1$ and $\l_2$, with the different choices of $h_{ij}$.

\subsection{The case when $h_{ij}=f\delta_{ij}$}
The evolution equation reads
\[
\p_t \rho_{ij}=f \Delta \rho_{ij}+f_i \rho_{jl,l}-f_{j}\rho_{il,l}
\]
Hence we compute
\[
\begin{split}
&\p_t |\rho|^2=f(\Delta |\rho|^2-2|\nabla \rho|^2)+2f_{i}\rho_{ij}\rho_{jl,l}\\
&\p_t u= f(\Delta u-|\nabla \rho^{+}|^2+|\nabla \rho^{-}|^2)+f_{i}(*\rho)_{ij}\rho_{jl,l}\\
&\p_t |\rho^{+}|^2=f(\Delta |\rho^{+}|^2-2|\nabla \rho^{+}|^2)+2f_i \rho^{+}_{ij} \rho_{jl, l}\\
&\p_t |\rho^{-}|^2=f(\Delta |\rho^{-}|^2-2|\nabla \rho^{-}|^2)+2f_i \rho^{-}_{ij} \rho_{jl, l}\\
\end{split}
\]
The conformal factor $f$ can take the form of $u^{-r}$ for some $r>0$, or something like $u^{-1}|\rho|^2$. The structure of the evolution equation along the flow indicates the almost identical technical difficulty to preserve the uniform positivity of $u$ (or $\l_2$) since the maximum principle argument fails similar as in the linear Hodge flow. It is still possible that the flow exists for all time with $f=u^{-r}$, but it is hard to obtain the uniform positivity of $u$ along the flow via a direct maximum principle argument. 

\subsection{The case when $h_{ij}=u^{-1} a_{ij}$} In this section we choose $g_\rho=u^{-1}g(A\cdot, A\cdot)$, which is the metric used in the Donaldson geometric flow. 
Correspondingly we have
\[
Q_{ij}=h_{ij}=\frac{a_{ij}}{u}. 
\]
Recall $a_{ij}=\rho_{ip}\rho_{jp}, b_{ij}=(*\rho)_{ip}(*\rho)_{jp},$ and  $a_{ij}+b_{ij}=|\rho|^2 \delta_{ij}$. We compute
\begin{equation}
\begin{split}
&\p_t |\rho|^2=\left(\frac{|\rho|^2 \rho_{kj}}{u}-(*\rho)_{kj}\right)\Delta \rho_{kj}+2\rho_{ij} h_{ik, j} \rho_{kl, l}\\
&\p_t u=\frac{1}{2} \rho_{kj}\Delta \rho_{kj} +(*\rho)_{ij} h_{ik, j} \rho_{kl,l}
\end{split}
\end{equation}
We compute, with $x=\frac{|\rho^{-}|^2}{|\rho^{+}|^2}$
\begin{equation}
\begin{split}
\p_t |\rho^{+}|^2=&\; \frac{1+x}{1-x} \Delta |\rho^{+}|^2+\frac{2}{1-x}\Delta |\rho^{-}|^2-\frac{2(1+x)}{1-x}|\nabla \rho^{+}|^2-\frac{4}{1-x} |\nabla \rho^{-}|^2\\
&+\left(\frac{b_{ik}}{u} +\delta_{ik}\right) \rho_{ij, j}\rho_{kl,l}+\left(\frac{|\rho|^2}{u}\right)_j \rho_{kj}\rho_{kl,l}\\
\p_t |\rho^{-}|^2=&\; \frac{2x}{1-x} \Delta |\rho^{+}|^2+\frac{1+x}{1-x} \Delta |\rho^{-}|^2-\frac{4x}{1-x}|\nabla \rho^{+}|^2-\frac{2(1+x)}{1-x} |\nabla \rho^{-}|^2\\
&+\left(\frac{b_{ik}}{u} -\delta_{ik}\right) \rho_{ij, j}\rho_{kl,l}+\left(\frac{|\rho|^2}{u}\right)_j \rho_{kj}\rho_{kl,l}\
\end{split}
\end{equation}
We compute, at points where  $|\rho^{-}|\neq 0$, 
\begin{equation}\label{a-1}
\begin{split}
\p_t \l_1=&\frac{\l_1}{\l_2}\Delta \l_1-\frac{\l_1}{\l_2}\left(J+K\right)+\frac{\l_1 b_{ik}-u\l_2 \delta_{ik}}{u(\l_1^2-\l_2^2)} \rho_{ij, j}\rho_{kl,l}+\left(\frac{\l_1^2+\l_2^2}{\l_1\l_2}\right)_j \frac{\l_1\rho_{kj}\rho_{kl,l}}{\l_1^2-\l_2^2}\\
\p_t \l_2=&\frac{\l_2}{\l_1}\Delta \l_2+\frac{\l_2}{\l_1}\left(K-J\right)+\frac{\l_1u\delta_{ik}-\l_2 b_{ik}}{u(\l_1^2-\l_2^2)} \rho_{ij, j}\rho_{kl,l}-\left(\frac{\l_1^2+\l_2^2}{\l_1\l_2}\right)_j \frac{\l_2\rho_{kj}\rho_{kl,l}}{\l_1^2-\l_2^2}
\end{split}
\end{equation}
We analyze \eqref{a-1} at points where $\l_1$ or $\l_2$ obtains its maximum or minimum respectively. Let $\l_1(p, t)=\max_{x\in M} \l_1(x, t)$. Then at the point $p$, we have
\[
\nabla \l_1=0, \nabla^2 \l_1\leq 0.
\]
This implies that $\nabla|\rho^+|+\nabla |\rho^{-}|=0$. We compute
\[
\frac{\langle \rho^{+}, \nabla \rho^{+}\rangle}{|\rho^{+}|}+\frac{\langle\rho^{-}, \nabla \rho^{-}\rangle}{|\rho^{-}|}=0
\]
We can choose a coordinate at $p$ such that $\rho_{12}=\l_1$ and $\rho_{34}=\l_2$ and \[\rho_{13}=\rho_{24}=\rho_{23}=\rho_{14}=0\]
This implies that $\nabla \rho_{12}=\nabla \l_1=0$, $\nabla \l_2=\nabla \rho_{34}$ at $p$. 
We compute
\[
J+K=\frac{\l_1}{\l_1^2-\l_2^2} (|\nabla \rho_{13}|^2+|\nabla \rho_{14}|^2+|\nabla \rho_{23}|^2+|\nabla \rho_{24}|^2)+\frac{2\l_2}{\l_1^2-\l_2^2} (\nabla \rho_{13}\nabla\rho_{24}-\nabla \rho_{23}\nabla \rho_{14})
\]
While for 
\[
\begin{split}
&\left(\frac{\l_1^2+\l_2^2}{\l_1\l_2}\right)_j \frac{\l_1\rho_{kj}\rho_{kl,l}}{\l_1^2-\l_2^2}=\frac{\l_1}{\l_2^2} (\rho_{34, 1}\rho_{2l, l}-\rho_{34, 2}\rho_{1l, l})+\frac{\rho_{34,3}\rho_{4l, l}-\rho_{34, 4}\rho_{3l, l}}{\l_2}\\
&\frac{\l_1 b_{ik}-u\l_2 \delta_{ik}}{u(\l_1^2-\l_2^2)} \rho_{ij, j}\rho_{kl,l}=\frac{1}{\l_2} [(\rho_{3j, j})^2+(\rho_{4j, j})^2]
\end{split}\]
It seems to be hard to bound $\l_1$ via a direct maximum principle argument since the term $\frac{\l_1}{\l_2^2} (\rho_{34, 1}\rho_{2l, l}-\rho_{34, 2}\rho_{1l, l})$ has the highest order when $\l_2<<1$, and this term does not have the favorite sign. 
We denote
\[
\p_t \l_2=\frac{\l_2}{\l_1}\Delta \l_2+\cR
\]
and analyze the behavior at minimum of $\l_2$, where $\l_2(q, t)=\min \l_2(\cdot, t)$. At the point $q$, where we choose a coordinate such that $\rho_{12}=\l_1$ and $\rho_{34}=\l_2$ and \[\rho_{13}=\rho_{24}=\rho_{23}=\rho_{14}=0\]
At $q$, 
$\nabla \l_2=0=\nabla \rho_{34}, \nabla \l_1=\nabla \rho_{12}$. We are not  able to prove that $\cR\geq 0$ at $q$.

\subsection{The case when $h_{ij}=u^{-2} a_{ij}$}

We have
\[
\rho_{ij} h_{ik}=\frac{|\rho|^2\delta_{kj}}{u^2}-\frac{(*\rho)_{kj}}{u}, (*\rho)_{ij}h_{ik}=\frac{\rho_{kj}}{u}
\]
We compute
\begin{equation}
\begin{split}
&\p_t |\rho|^2=\left(\frac{|\rho|^2\rho_{kj}}{u^2}-\frac{(*\rho)_{kj}}{u} \right) \Delta \rho_{kj}+2\left(\frac{|\rho|^2\rho_{kj}}{u^2}-\frac{(*\rho)_{kj}}{u}\right)_j \rho_{kl, l}-2h_{ik} \rho_{ij, j}\rho_{kl,l}\\
&\p_t u=\frac{\rho_{kj}}{2u} \Delta \rho_{kj}+\left(\frac{\rho_{kj}}{u}\right)_j \rho_{kl, l} 
\end{split}
\end{equation}
We compute
\begin{equation*}
\begin{split}
\p_t (|\rho|^2+2u)=&\left(\frac{|\rho|^2+u}{u^2} \rho_{kj}-\frac{(*\rho)_{kj}}{u}\right)\Delta \rho_{kj}+2\left(\frac{|\rho|^2+u}{u^2} \rho_{kj}-\frac{(*\rho)_{kj}}{u}\right)_j \rho_{kl,l}-2h_{ik} \rho_{ij, j}\rho_{kl,l}\\
\p_t (|\rho|^2-2u)=&\left(\frac{|\rho|^2-u}{u^2} \rho_{kj}-\frac{(*\rho)_{kj}}{u}\right)\Delta \rho_{kj}+2\left(\frac{|\rho|^2-u}{u^2} \rho_{kj}-\frac{(*\rho)_{kj}}{u}\right)_j \rho_{kl,l}-2h_{ik} \rho_{ij, j}\rho_{kl,l}
\end{split}
\end{equation*}
We compute
\begin{equation}
\begin{split}
\p_t |\rho^{+}|^2=&\frac{|\rho|^2}{2u^2}\Delta |\rho^{+}|^2+\frac{|\rho|^2+2u}{2u^2}\Delta |\rho^{-}|^2-\frac{|\rho|^2}{u^2} |\nabla \rho^{+}|^2-\frac{|\rho|^2+2u}{u^2}|\nabla \rho^{-}|^2\\
&+\left(\frac{|\rho|^2+u}{u^2} \rho_{kj}-\frac{(*\rho)_{kj}}{u}\right)_j \rho_{kl,l}-h_{ik} \rho_{ij, j}\rho_{kl,l}\\
\p_t |\rho^{-}|^2=&\frac{|\rho|^2-2u}{2u^2}\Delta |\rho^{+}|^2+\frac{|\rho|^2}{2u^2}\Delta |\rho^{-}|^2-\frac{|\rho|^2-2u}{u^2} |\nabla \rho^{+}|^2-\frac{|\rho|^2}{u^2}|\nabla \rho^{-}|^2\\
&+\left(\frac{|\rho|^2-u}{u^2} \rho_{kj}-\frac{(*\rho)_{kj}}{u}\right)_j \rho_{kl,l}-h_{ik} \rho_{ij, j}\rho_{kl,l}\\
\end{split}
\end{equation}
It follows that
\begin{equation}
\begin{split}
\p_t \l_1=&\frac{\Delta \l_1}{\l_2^2}-\frac{1}{\l_2^2} (J+K)+\frac{\l_1 b_{ik}-u\l_2 \delta_{ik}}{u^2(\l_1^2-\l_2^2)} \rho_{ij, j}\rho_{kl,l}+P_1\\
\p_t \l_2=&\frac{\Delta \l_2}{\l_1^2}-\frac{1}{\l_1^2} (K-J)+\frac{\l_1u\delta_{ik}-\l_2 b_{ik}}{u^2(\l_1^2-\l_2^2)} \rho_{ij, j}\rho_{kl,l}+P_2
\end{split}
\end{equation}
where we denote
\[
\begin{split}
P_1=&\frac{\l_{1, j}\rho_{kl, l}}{\l_1(\l_1^2-\l_2^2)}\left(\frac{(*\rho)_{kj}}{\l_2}-\frac{\rho_{kj}}{\l_1}\right)+\frac{\l_{2, j}\rho_{kl, l}}{\l_2(\l_1^2-\l_2^2)}\left(\frac{(*\rho)_{kj}}{\l_2}+\frac{\rho_{kj}}{\l_1}-\frac{2\l_1}{\l_2^2}\rho_{kj}\right),\\
P_2=&\frac{\l_{1, j}\rho_{kl, l}}{\l_1(\l_1^2-\l_2^2)}\left(\frac{(*\rho)_{kj}}{\l_2}+\frac{\rho_{kj}}{\l_2}-\frac{2\l_2}{\l_1^2}\rho_{kj}\right)+\frac{\l_{2, j}\rho_{kl, l}}{\l_2(\l_1^2-\l_2^2)}\left(\frac{(*\rho)_{kj}}{\l_2}-\frac{\rho_{kj}}{\l_1}\right)
\end{split}
\]

\subsection{The case when $h_{ij}=u^{-1}b_{ij}$}

In this section we choose
\[
Q_{ij}=h_{ij}=\frac{b_{ij}}{u}. 
\]
We have
\[
\rho_{ij} h_{ik}=(*\rho)_{kj},\; (*\rho)_{ij} h_{ik}=\frac{|\rho|^2(*\rho)_{kj}-u\rho_{kj}}{u} .
\]
We compute
\begin{equation}
\begin{split}
&\p_t |\rho|^2=(*\rho)_{kj}\Delta \rho_{kj}-\frac{2b_{ik}}{u}\rho_{ij, j}\rho_{kl,l}\\
&\p_t u=\frac{|\rho|^2(*\rho)_{kj}-u\rho_{kj}}{2u} \Delta \rho_{kj}+\left(\frac{|\rho|^2(*\rho)_{kj}-u\rho_{kj}}{u}\right)_j \rho_{kl, l}
\end{split}
\end{equation}
Hence we have
\begin{equation}
\begin{split}
\p_t |\rho^{+}|^2=&\frac{(|\rho|^2+u)(*\rho)_{kj}-u\rho_{kj}}{2u} \Delta \rho_{kj}-\frac{b_{ik}}{u}\rho_{ij, j}\rho_{kl,l}+\left(\frac{|\rho|^2(*\rho)_{kj}-u\rho_{kj}}{u}\right)_j \rho_{kl, l}\\
\p_t |\rho^{-}|^2=&\frac{(u-|\rho|^2)(*\rho)_{kj}+u\rho_{kj}}{2u} \Delta \rho_{kj}-\frac{b_{ik}}{u}\rho_{ij, j}\rho_{kl,l}-\left(\frac{|\rho|^2(*\rho)_{kj}-u\rho_{kj}}{u}\right)_j \rho_{kl, l}\\
\end{split}
\end{equation}
Hence we have
\begin{equation}\label{L0-1}
\begin{split}
\p_t \l_1=\frac{\l_2}{\l_1}\Delta \l_1-\frac{\l_2}{\l_1} (J+K)-\frac{\l_1 b_{ik}-u\l_2 \delta_{ik}}{u(\l_1^2-\l_2^2)} \rho_{ij, j}\rho_{kl,l}-R_1\\
\p_t \l_2=\frac{\l_1}{\l_2}\Delta \l_2-\frac{\l_1}{\l_2} (K-J)+\frac{\l_2 b_{ik}-u\l_1 \delta_{ik}}{u(\l_1^2-\l_2^2)} \rho_{ij, j}\rho_{kl,l}+R_2
\end{split}
\end{equation}
where we denote
\[
\begin{split}
R_1=\left(\frac{\l_1^2+\l_2^2}{\l_1\l_2}\right)_j \frac{\l_2 (*\rho)_{kj}\rho_{kl, l}}{\l_1^2-\l_2^2}, \; R_2=\left(\frac{\l_1^2+\l_2^2}{\l_1\l_2}\right)_j \frac{\l_1 (*\rho)_{kj}\rho_{kl, l}}{\l_1^2-\l_2^2}
\end{split}
\]

\subsection{The case when $h_{ij}=u^{-2}b_{ij}$}

 In this section we choose
\begin{equation}
Q_{ij}=b_{ij}, h_{ij}=\frac{b_{ij}}{u^2}
\end{equation}
We compute,
\[
\rho_{ij} h_{ik}=\frac{(*\rho)_{kj}}{u},\; (*\rho)_{ij} h_{ik}=\frac{|\rho|^2(*\rho)_{kj}-u\rho_{kj}}{u^2} .
\]
Hence we have the following, 
\begin{equation}
\begin{split}
&\p_t |\rho|^2=\frac{(*\rho)_{kj}}{u} \Delta \rho_{kj}+2 \left(\frac{(*\rho)_{kj}}{u} \right)_{,\; j}\rho_{kl, l}-2 h_{ik}\rho_{ip, p}\rho_{kl, l}\\
&\p_t (2u)=\frac{|\rho|^2(*\rho)_{kj}-u\rho_{kj}}{u^2} \Delta \rho_{kj}+2\left(\frac{|\rho|^2(*\rho)_{kj}-u\rho_{kj}}{u^2} \right)_{,\; j}\rho_{kl, l}\\
&\p_t (|\rho|^2+2u)=R_{kj}\Delta \rho_{kj}+2R_{kj,\,j} \rho_{kl, l}-2h_{ik}\rho_{ip, p}\rho_{kl, l}, \\
&\p_t (|\rho|^2-2u)=S_{kj}\Delta \rho_{kj}+2S_{kj,\;j} \rho_{kl, l}-2h_{ik}\rho_{ip, p}\rho_{kl, l}\\
\end{split}
\end{equation}
where we denote
\[
R_{kj}=\frac{(|\rho|^2+u)(*\rho)_{kj} -u\rho_{kj}}{ u^2},\;\;S_{kj}=\frac{\rho_{kj}+(*\rho)_{kj}}{u}-\frac{|\rho|^2}{u^2} (*\rho)_{kj}
\]
Recall that
\[
\l_1= \frac{\sqrt{|\rho|^2+2u}}{2}+\frac{\sqrt{|\rho|^2-2u}}{2}, \l_2=\frac{\sqrt{|\rho|^2+2u}}{2}-\frac{\sqrt{|\rho|^2-2u}}{2}
\]
Note that $\l_1$ and $\l_2$ are both smooth functions  when $|\rho^{-}|^2\neq 0$. 
We compute
\begin{equation}\label{L1-1}
\begin{split}
\p_t \l_1=&\frac{\p_t (|\rho|^2+2u)}{4\sqrt{|\rho|^2+2u}}+\frac{\p_t (|\rho|^2-2u)}{4\sqrt{|\rho|^2-2u}}\\
=& \left(\frac{R_{kj}\Delta \rho_{kj}}{4\sqrt{|\rho|^2+2u}}+\frac{S_{kj}\Delta \rho_{kj}}{4\sqrt{|\rho|^2-2u}}\right) + \left(\frac{R_{kj, j}\rho_{kl, l}}{2\sqrt{|\rho|^2+2u}}+\frac{S_{kj, j}\rho_{kl, l}}{2\sqrt{|\rho|^2-2u}}\right) \\
&-\left(\frac{1}{2\sqrt{|\rho|^2+2u}}+\frac{1}{2\sqrt{|\rho|^2-2u}}\right) h_{ik} \rho_{ip, p}\rho_{kl,l}
\end{split}
\end{equation}
and similarly 
\begin{equation}\label{L2-1}
\begin{split}
\p_t \l_2=&\frac{\p_t (|\rho|^2+2u)}{4\sqrt{|\rho|^2+2u}}-\frac{\p_t (|\rho|^2-2u)}{4\sqrt{|\rho|^2-2u}}\\
=& \left(\frac{R_{kj}\Delta \rho_{kj}}{4\sqrt{|\rho|^2+2u}}-\frac{S_{kj}\Delta \rho_{kj}}{4\sqrt{|\rho|^2-2u}}\right) + \left(\frac{R_{kj, j}\rho_{kl, l}}{2\sqrt{|\rho|^2+2u}}-\frac{S_{kj, j}\rho_{kl, l}}{2\sqrt{|\rho|^2-2u}}\right) \\
&-\left(\frac{1}{2\sqrt{|\rho|^2+2u}}-\frac{1}{2\sqrt{|\rho|^2-2u}}\right) h_{ik} \rho_{ip, p}\rho_{kl,l}
\end{split}
\end{equation}
By \eqref{L1-1}, we compute, using $\rho=\rho^{+}+\rho^{-}$ and $*\rho=\rho^{+}-\rho^{-}$, 
\begin{equation}\label{L1-2}
\begin{split}
\p_t \l_1=&\frac{\Delta \l_1}{\l_1^2}-\frac{1}{\l_1^2}(J+K)-\frac{\l_1  h_{ik}\rho_{ip,p}\rho_{kl,l}}{\l_1^2-\l_2^2}+\frac{|\rho_{kl,l}|^2}{\l_1(\l_1^2-\l_2^2)}\\
&+\frac{\l_{1, j}\rho_{kl, l}}{\l_1 (\l_1^2-\l_2^2)}\left( \left(\frac{2\l_2}{\l_1^2}-\frac{1}{\l_2}\right)(*\rho)_{kj}-\frac{\rho_{kj}}{\l_1}\right)+\frac{\l_{2, j} \rho_{kl, l}\left(\frac{(*\rho)_{kj}}{\l_2}-\frac{\rho_{kj}}{\l_1}\right)}{\l_2(\l_1^2-\l_2^2)}
\end{split}
\end{equation}
and similarly we compute
\begin{equation}
\label{L2-2}
\begin{split}
\p_t \l_2=&\frac{\Delta \l_2}{\l_2^2}+\frac{1}{\l_2^2}(K-J)+\frac{\l_2  h_{ik}\rho_{ip,p}\rho_{kl,l}}{\l_1^2-\l_2^2}-\frac{|\rho_{kl,l}|^2}{\l_2(\l_1^2-\l_2^2)}\\
&+\frac{\l_{2, j}\rho_{kl, l}}{\l_2 (\l_1^2-\l_2^2)}\left(\frac{\rho_{kj}}{\l_2}+\frac{(*\rho)_{kj}}{\l_1}-\frac{2\l_1(*\rho)_{kj}}{\l_2^2}\right)+\frac{\l_{1, j} \rho_{kl, l}\left(\frac{\rho_{kj}}{\l_2}-\frac{(*\rho)_{kj}}{\l_1}\right)}{\l_1(\l_1^2-\l_2^2)}
\end{split}
\end{equation}

When $\rho=\l_1 (x_1, x_2) dx_1\wedge dx_2+\l_2(x_3, x_4) dx_3\wedge dx_4$, we have already seen that the nonlinear Hodge flow with $h_{ij}=u^{-2} b_{ij}$ preserves the structure and satisfies
\begin{equation}\label{L3}
\p_t \l_1=\Delta_{x_1, x_2} (-\l_1^{-1}), \p_t \l_2=\Delta_{x_3, x_4}(-\l_2^{-1})
\end{equation}
It is a direct computation to check that \eqref{L3} coincides with \eqref{L1-2} and \eqref{L2-2}. For example, $J=K=0$ and $\rho_{1l, l}=\p_{x_2}\l_1, \rho_{2l, l}=-\p_{x_1}\l_1, \rho_{3l, l}=\p_{x_4}\l_2, \rho_{4l, l}=-\p_{x_3}\l_2$. 
Note that \eqref{L3} is a fast diffusion equation for $\l_1$ and $\l_2$ respectively.  The maximum principle indeed applies to show that $\l_1, \l_2$ remain uniformly bounded and uniformly positive. By the parabolic theory the flow exists for all time with uniform control, and it will converge to $\omega=dx_1\wedge dx_2+dx_3\wedge dx_4$. However, the maximum principle argument does not apply directly to \eqref{L1-2} and \eqref{L2-2} for the general case. 
On the other hand, \eqref{L2-2} might still suggest that $\l_2$  stays uniformly positive along the flow. When $\l_2$ approaches to zero, the leading term $\l_2^{-2}\Delta \l_2$ is the same as fast diffusion equation and it might force $\l_2$ to stay uniformly positive. Nevertheless it is clear that the system is really complicated and it remains a substantial challenge to address these difficulties.  

\bibliographystyle{plain}

\end{document}